\documentclass[12pt]{article}
\usepackage{setspace}
\usepackage{graphicx} 
\usepackage{todonotes}
\usepackage{amsmath,amssymb}
\usepackage[a4paper, margin=2cm]{geometry}
\usepackage{booktabs}
\usepackage[inkscapelatex=false]{svg}
\usepackage{verbatim}
\usepackage{makecell}
\usepackage{soul} 
\usepackage[english]{babel}
\renewcommand{\todo}[1]{}
\renewcommand{\todo}[2][]{}
\usepackage{hyperref} 
\usepackage{csquotes}
\usepackage{lmodern}
\usepackage{tikz}
\usetikzlibrary{shapes.geometric, arrows}
\usetikzlibrary{calc}

\usepackage[
  backend=biber,
  language=english,
  style=vancouver,   
  sorting=none,      
  giveninits=true,   
  doi=true,          
  url=false,         
  eprint=false,
  isbn=false,
  maxnames=10,       
  minnames=1         
]{biblatex}
\DeclareNameAlias{author}{family-given}
\DeclareFieldFormat{year}{#1}
\DeclareFieldFormat{date}{%
  \printfield{year}%
}
\DeclareFieldFormat[article]{volume}{\mkbibbold{#1}}
\DeclareFieldFormat{volume}{\mkbibbold{#1}}
\DeclareFieldFormat[inproceedings]{booktitle}{\mkbibemph{#1}}
\DeclareFieldFormat{eventtitle}{\mkbibemph{#1}}
\DeclareFieldFormat[book]{title}{\mkbibemph{#1}}
\DeclareFieldFormat[incollection]{title}{\mkbibemph{#1}}
\DeclareFieldFormat[inbook]{title}{\mkbibemph{#1}}

\DeclareFieldFormat{labelnumberwidth}{\mkbibbrackets{#1}}
\usepackage{yhmath}
\DeclareMathOperator{\diag}{diag}

\usepackage{multicol}

\begin{document}

\title{Markov chain Monte Carlo for Bayesian inference of the non-conducting region in intra-atrial reentrant tachycardia}

\author{Maarten Volkaerts\footnote{KU Leuven, Department of Computer Science (NUMA), Celestijnenlaan 200A, 3001 Leuven, Belgium,
		maarten.volkaerts@kuleuven.be}, 
	Marie Cloet$^\ast$, 
	Piet Claus\footnote{KU Leuven, Department of Cardiovascular Sciences (Cardiovascular Imaging and Dynamics), UZ Herestraat 49, 3000 Leuven, Belgium,}, Hans Dierckx\footnote{Leiden University Medical Center, Department of Cardiology}, Giovanni Samaey$^\ast$}
\maketitle
\begin{abstract}
We present a Bayesian approach to estimate the parameters of mathematical models of cardiac electrophysiology with quantified uncertainty. Such models capture the dynamics of the electrical signal that coordinates the muscle cell contraction in the heart wall and can support cardiac arrhythmia treatment. We consider an illustrative case motivated by a cardiac arrhythmia, namely, by intra-atrial reentrant tachycardia. We estimate a low-dimensional geometrical parameter that describes the boundary of an electrically non-conducting region in the heart tissue from synthetic electrical measurements outside of the tissue. Instead of relying on a deterministic fit for this region, we estimate a posterior distribution on the geometrical parameter using Bayesian inference that captures the uncertainty due to measurement errors. We propose a likelihood based on a set of quantities that characterize the data for improved accuracy. To efficiently approximate the posterior distribution, we propose a compressed likelihood function and an adapted Metropolis-Hastings (MH) algorithm. We obtain an algorithm that strongly decreases the number of samples by using an adaptive proposal strategy. Our algorithm also gives attention to the impact of discretization errors on inference outcomes, as these introduce artificial discontinuities in the posterior if not properly addressed. We account for discretization errors in the likelihood and in the accept-reject step of our adapted MH algorithm to improve the robustness of our estimates and to further increase the sampling efficiency. All of these elements combined give us a method that efficiently estimates the non-conducting parameters with uncertainty. We perform several experiments with different amounts of measurement noise and illustrate how this translates into the posterior distributions.
\end{abstract}
\section{Introduction}
Personalized heart models, also known as virtual twins of the heart, are gaining increasing attention due to their potential to enhance understanding and prediction of cardiac function \cite{niederer2019,niederer2016,trayanova2018}. By integrating biophysical principles with data, they offer a virtual representation of cardiac electrophysiology that can improve treatment strategies for cardiac arrhythmia \cite{arevalo2016,ashikaga2013,boyle2019,deng2016,mcdowell2015,prakosa2018,relan2011}.

The calibration of these models involves anatomical and functional components, describing the anatomy and the electrophysiological and mechanical properties of the heart tissue, respectively. Due to the limited and noisy data, assumptions and approximations are made along the calibration pipeline, for example, in assigning averaged tissue properties or in image segmentation \cite{niederer2019}. 

A key challenge is to understand how such assumptions and data limitations affect model predictions. This has motivated research efforts on uncertainty quantification (UQ) during model calibration and prediction \cite{mirams2016}. The objective is to associate model outputs with quantitative measures of reliability, allowing for the systematic assessment of prediction quality. Such measures are essential not only for potential clinical use but also for developing a trustworthy computational modeling pipeline.

\subsection{Bayesian estimation of non-conducting regions}
This work develops a Bayesian methodology for uncertainty-aware estimation of a geometrical parameter that characterizes an electrically non‑conducting region in the cardiac tissue from synthetic electrogram (EGM) data. We consider this case as an introductory case for our Bayesian methodology for parameter estimation with uncertainty in mathematical models of cardiac electrophysiology (EP). 

Mathematical models of cardiac EP describe how electrical waves propagate through heart tissue, triggering the activation of muscle cells and thereby coordinating cardiac contraction. Because EP dynamics are tightly linked to the mechanisms underlying cardiac arrhythmias, reliable personalization of these models is essential for constructing patient‑specific cardiac models. Such personalization requires identifying the electrical properties of the tissue. However, this task is challenging due to the large number of parameters involved and their spatial (extrinsic) and temporal (intrinsic) variability \cite{clayton2020}. Addressing this difficulty by including uncertainties in the resulting models is exactly what we aim to achieve. 

We carefully drafted a setup to obtain an illustrative case for the introduction of our Bayesian methodology that can be extended to more complicated settings in follow-up work \cite{volkaerts2024}. Our long-term goal is achieve feasible inference of high-dimensional parameters in detailed EP models. As these models are based on time-dependent partial differential equations (PDEs), we employ the simplest PDE-based model. We partition the tissue into two categories: healthy and non‑conducting. This offers a simplified representation of spatial heterogeneity, which can still be linked to a clinically relevant arrhythmia (see Section \ref{IART}). Our aim is to estimate a geometrical parameter that defines the boundary between these two regions. This formulation yields a low‑dimensional parameter estimation problem, which can be extended to high‑dimensional parameters in subsequent research. 

While classical optimization seeks an optimal fit only, Bayesian inference defines a posterior probability distribution on the geometrical parameter, providing both an estimate and a measure of the uncertainty associated with this estimate. More specifically, Bayesian inference models the geometrical parameter $\vartheta$, the model output $y$ and the discrepancy $\varepsilon$ due to measurement error as random variables
\begin{equation} \label{eq:misfit}
y = h(\mathcal{F}(\vartheta))+\varepsilon,
\end{equation}
where $\mathcal{F(\vartheta)}$ typically involves solving the model for the parameter value $\vartheta$ and $h(.)$ is the observation map that projects this solution to what is observed in the data.

Using Bayes' theorem, we can identify the posterior probability distribution on $\vartheta$ that captures the resulting uncertainties on the geometrical parameter as follows \cite{law2015}:
\begin{equation}\label{eq:bayesianinference}
P(\vartheta|y)= \frac{\mathcal{L}(y|\vartheta)\pi(\vartheta)}{P(y)},
\end{equation} 
where $\pi(\vartheta)$ is the prior distribution encoding prior knowledge about the parameter, $\mathcal{L}(y|\vartheta)$ is the likelihood quantifying the probability of the discrepancy between model output for $\vartheta$ and observed data, and $P(y)$ is the evidence, an unknown normalization constant. 
\subsection{Motivational case: Intra-Atrial Reentrant Tachycardia} \label{IART}
We employ a mathematical model of cardiac EP that models intra-atrial reentrant tachycardia (IART). IART is a cardiac arrhythmia characterized by an increased frequency of electrical activation of the muscle cells in the atrium due to a recurring electrical circuit within this atrium \cite{homoud2025}. While many recurring electrical circuits, also called reentries, can result in such an arrhythmia, we specifically focus on the case where this reentry occurs around a scarred and, consequently, non-conducting region.

\begin{figure}[ht]
    \centering
    \includegraphics[width=0.45\linewidth]{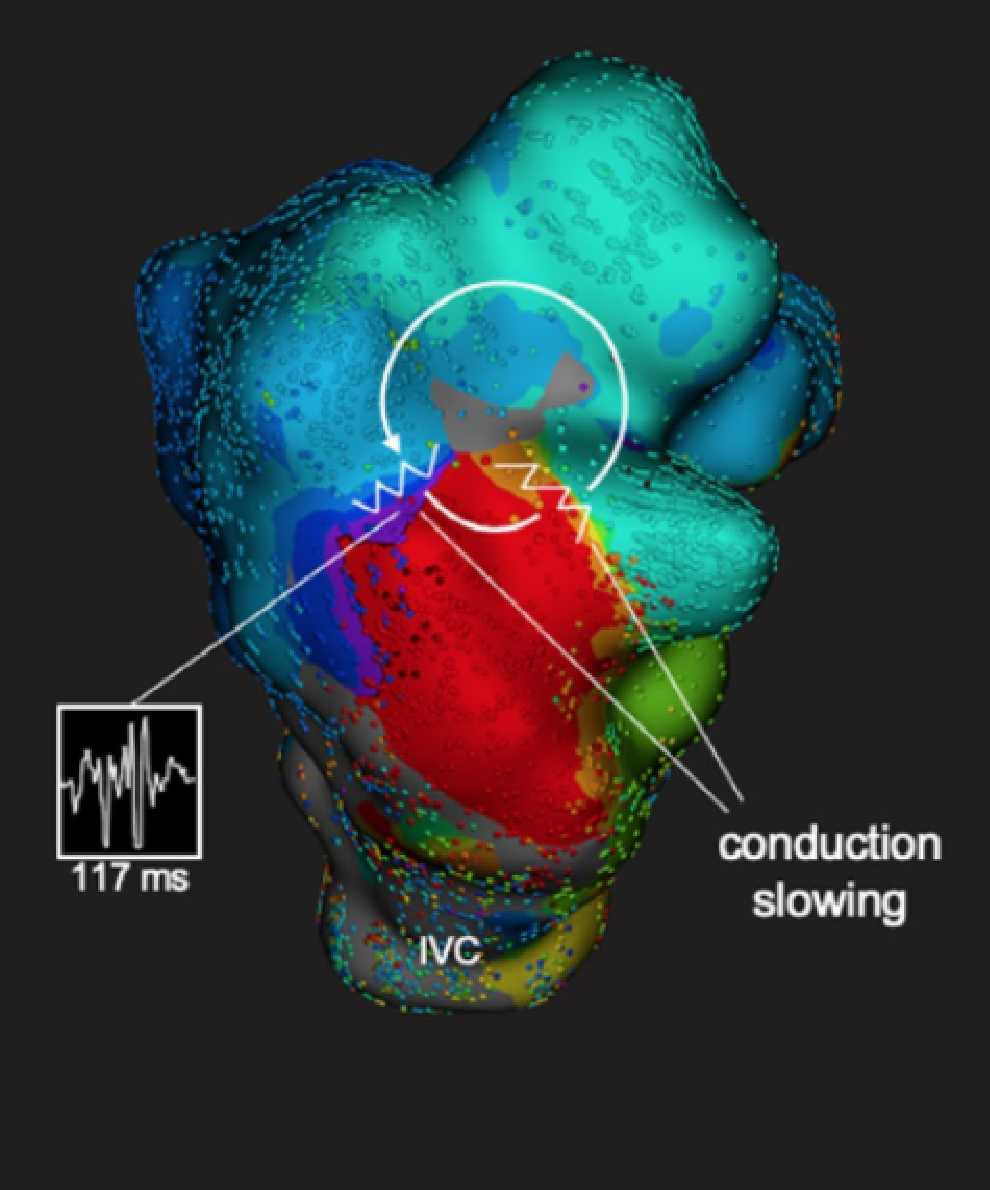}
    \caption{The colors annotate the arrival time of the electrical wave from early (red) to late (dark blue) under IART. The grey area is a surgical scar and does not yield electrical activity. We observe reentry around this scar. Adapted from \cite{kahle2021}. IVC = inferior vena cava.}
    \label{fig:IART}
\end{figure}

Figure \ref{fig:IART} shows the electrical propagation during such a reentry. The electrical wave moves along the color gradient and short-circuits around this non-conducting region, shown in gray. Unlike during a healthy heart rhythm, where a planar electrical wave moves through the atria, activating them once per heartbeat, this reentry causes continuous activation of parts of the atrium. 

In this paper, we compose a simplified mathematical model to describe cardiac EP under an IART condition.  IART occurring around a non-conducting region serves as an ideal first case study due to its clinical relevance and moderate complexity.  The spatial heterogeneity of the affected tissue is limited, allowing for characterization by a low-dimensional parameter. Moreover, the arrhythmia dynamics are relatively straightforward and amenable to mathematical modeling. 

\subsection{Related work}
In this section, we provide a concise overview of related work. Given the multidisciplinary nature of this project, it intersects with several research domains. We concentrate on the areas most relevant to the methodological contributions of this study: the calibration of EP models and the application of Markov chain Monte Carlo methods.
\subsubsection{Calibration of EP models}
Calibration of mathematical models of EP is often limited to a personalized division of the tissue into healthy and unhealthy regions, each of which is imposed with averaged electrical properties \cite{arevalo2016,boyle2019,prakosa2018}, or estimation of a low-dimensional description of the electrical properties \cite{gillette2021}. Our work aligns with the former approach, but could, without additional complexity, be applied to the latter as well. We improve upon earlier work on personalized tissue division by additionally quantifying the uncertainty in this division, arising from limited data availability and other discrepancies, using Bayesian inference.

The interest in UQ for personalized heart models has resulted in several other studies that perform UQ during the estimation of the parameters of mathematical models of EP. A common denominator is that they all alleviate the sampling cost by replacing the forward model with an emulator. A first study employed this to infer a space-dependent conductivity \cite{konukoglu2011}. Recently, a collection of studies paved the way to enabling estimation of the parameters of a mathematical model of EP in the atrium by using Gaussian processes as surrogates for the model to alleviate sampling cost \cite{coveney2020a,coveney2020,coveney2021,coveney2022}. Our work differentiates from these studies with its focus on using the original PDE-based models instead of surrogates to make inference using these models feasible. Another application of UQ that closely matches our work is a forward UQ study that evaluated how the reentry around a non-conducting region is influenced by the uncertainty on the geometry of this region \cite{clayton2024}. 

Finally, our method potentially improves upon the current clinical estimation of non-conducting regions. Bayesian inference postulates a rigorous mathematical framework that integrates the forward model encapsulating the dynamics, knowledge of the measurement noise and the data to enable a consistent and interpretable link between this region and the surrounding reentry. This contrasts with current practice, where catheters are manually navigated across the endocardial surface to construct activation maps and identify low-voltage regions as non-conducting. Such mappings do not guarantee consistency between the inferred non-conducting region and the observed activation patterns.
\subsubsection{Markov chain Monte Carlo methods}

Full characterization of the posterior distribution in Bayesian inference problems is typically achieved using a Markov chain Monte Carlo (MCMC) method. The main challenge of this method is its high computational cost. Estimation of the posterior necessitates a large number of samples, and each sample requires solving the forward model to evaluate the likelihood function. Typically, this forward model is an expensive partial differential equation (PDE) such as the monodomain model of cardiac EP in our study. This has resulted in a vast literature that aims to alleviate this cost by reducing the number of required samples or the cost per likelihood evaluation. 

The cost of the likelihood evaluation can be reduced by replacing the forward simulation with a faster, approximate solution. This has resulted in multifidelity \cite{peherstorfer2018b} and multilevel methods \cite{dodwell2019,lykkegaard2023,vanmechelen2025}. In this work, we consider the influence of the discretization error on the inference results, which relates closely to the multilevel strategies. Discretization error can seriously affect the results of the inference, as it causes an extra source of discrepancy between model output and data that is often ignored. A typical result of this effect is overconfidence in wrong predictions \cite{poot2025}. Recently, a review compared several approaches to model the discretization error during Bayesian inference for finite element method discretizations \cite{poot2025}. 

Decreasing the number of samples is achieved by choosing a proposal distribution that optimizes the parameter space exploration and minimizes the autocorrelation of the Markov chains. The classical random-walk Metropolis-Hastings (MH) faces important limitations, including poor sampling of distributions that are multimodal, high-dimensional or have an anisotropic covariance matrix. This has resulted in research into a variety of alternative proposal strategies \cite{chen2019b,cotter2013,cui2016,haario2006,haario2001}. The main challenge encountered in this work concerns the sampling of posterior distributions with an anisotropic covariance. This challenge can be resolved by using the Adaptive Metropolis method, an extension of the MH algorithm that adapts the covariance of the Gaussian proposal during sampling based on the covariance structure of the collected samples \cite{haario2001}. 

To minimize the number of samples, we propose an adapted MH algorithm that implements a proposal inspired by the AM method and an accept-reject step that can deal with discretization errors.

\subsection{Outline \& Contributions}
In this paper, we focus on fast and reliable UQ during the inference of a geometrical parameter, which describes a non-conducting region, in mathematical models of EP models from synthetic electrogram data. First, we discuss the exact inference problem in detail, describing the forward model and the data in Sections \ref{forward-model} and \ref{data}, respectively. In the remainder of the paper, we make the following contributions:
\begin{enumerate} 
    \item In Section \ref{Bayes}, we formulate the geometrical parameter estimation as a Bayesian inference problem to facilitate uncertainty quantification.
    \item In Section \ref{summary-statistic}, we reduce the dimension of the data from a large number of time traces to a low-dimensional vector of quantities that are characteristic of the underlying dynamics. Using a compressed likelihood based on these quantities improves its accuracy and enables the application of our methodology to clinical data in follow-up research.
    \item We establish the role of discretization error of the forward model on the results of Bayesian inference of geometrical parameters and account for it in the likelihood in Section \ref{disc-ll}.
    \item In Sections \ref{MCMC} and \ref{results-exp-1}, we outline and test our adapted MH method. We adapt the proposal and the accept-reject step of the classical MH algorithm to achieve efficient and robust sampling under parameter correlations and discretization error, respectively.
    \item Finally, we illustrate the benefits of stochastic parameter estimation compared to deterministic approaches in Sections \ref{results-exp-2} and \ref{results-exp-3}. In contrast to the latter offering one single fit, our method includes a reliability measure by quantifying the uncertainty. We demonstrate the added value in an example, where the structure of the uncertainty offers insights into identifiability issues.
\end{enumerate}
The code that implements the numerical methodology and performs the numerical experiments is published on gitlab \footnote{https://gitlab.kuleuven.be/numa/public/mcmc-for-bayesian-inference-of-non-conducting-region-in-iart}.


\section{Mathematical model of EP} \label{forward-model}
 As the forward model for Bayesian inference, we use the monodomain model. The goal is not to capture detailed physiology but to provide a controlled and computationally efficient setting for evaluating UQ performance during parameter inference. Section \ref{tissue} presents this model that captures the electrical wave propagation through cardiac tissue. Section \ref{IART-model} explains how this model can replicate IART scenarios by adding a non-conducting region and imposing an appropriate initial condition. In Section \ref{personalized-IART}, we define the geometrical parameter that models the non-conducting region and that is inferred from data in this work. Then, Section \ref{numerical-solution} describes the numerical methods used to solve the model. Finally,  in Section \ref{prepacing}, we elaborate on the prepacing simulation that generates the initial condition for our mathematical model.
\subsection{Tissue model}\label{tissue} \
The dynamics of the electrical wave are described using the monodomain model with the modified Mitchell-Schaeffer cell model \cite{corrado2016, mitchell2003}:

  \begin{align}
\label{eq:model}
        		\frac{\partial V_m(x,t)}{\partial t}&= \nabla D  \nabla V_m(x,t) + h(x,t) \frac{ V_m(x,t)(V_m(x,t)-v_{\text{gate}})(1-V_m(x,t))}{\tau_{\text{in}}}  \\
  & \qquad - (1-h(x,t))\frac{V_m(x,t)}{\tau_{\text{out}} }\nonumber\\
    		\label{eq:h}\frac{\partial h(x,t)}{\partial t} &= \begin{cases} \frac{1-h(x,t)}{\tau_{\text{open}}}, & \mbox{if } V_m(x,t) < v_{\text{gate}}\\ \frac{-h(x,t)}{\tau_{\text{close}}}, & \mbox{if }V_m(x,t)>v_{\text{gate}}\\ \end{cases}.
\end{align}
The normalized transmembrane potential, $V_m(x,t)$, models local electrical activation in the tissue. The monodomain model incorporates both a diffusion term, which governs the conduction of the electrical wave through the tissue, and a reaction term, which captures cellular EP. The modified Mitchell--Schaeffer model is a cell model consisting of a contribution from all inward ion currents and one for all outward currents, with a slow gating variable $h$ regulating the ion channels opening and closing. 

The monodomain model with Mitchell-Schaeffer cell dynamics captures the overall qualitative cardiac electrical behavior while minimizing the model complexity. As the model must be solved a large number of times to estimate the posterior distribution of the parameters (see Section \ref{MCMC}), its minimal complexity keeps the computational time feasible. Large research efforts have gone into UQ methods that reduce the computational time without compromising the accuracy of the inference results by replacing a large number of forward simulations with faster, approximate simulations. However, these methods are yet to be tailored to cardiac electrophysiology models. The methodology presented in this paper provides a foundation for integrating such techniques in follow-up research, thereby enabling the solution of the same inference problem using more detailed cardiac EP models within a feasible time frame. 

The time constants of the Mitchell-Schaeffer model are set to the following values: $\tau_{\text{in}}=0.3$~ms, $\tau_{\text{out}}=6.0$~ms, $\tau_{\text{open}}=120$~ms, and $\tau_{\text{close}}= 150$~ms. The intracellular and extracellular conductivity are both isotropic, with $\sigma_i = \sigma_{il}= \sigma_{it}=0.174$~mS/mm and $\sigma_e = \sigma_{el}= \sigma_{et}=0.625$~mS/mm. In the monodomain case, we set the overall conductivity to the harmonic mean of both conductivities: ${\sigma = \frac{\sigma_i\sigma_e}{\sigma_i + \sigma_e}}$. The diffusion constant in the monodomain model captures the conductivity with ${D_{\text{healthy}} = \frac{\sigma}{\chi C_m}}$ using standardized values for the surface-to-volume ratio of the cells $\chi=140$~$\text{mm}^{-1}$ and the membrane capacitance $C_m = 0.01$~$\mu$F/$\text{mm}^2$.
\subsection{Model of IART} \label{IART-model}
In this study, we aim to simulate intra-atrial reentrant tachycardia  (IART). To this end, we define the model from Section \ref{tissue} over a 2D tissue slab $\Omega_{\text{atrium}}=[0~\text{cm},10~\text{cm}]\times [0~\text{cm},10~\text{cm}]$, as an approximation for the atrium. Reentry occurs around a non-conducting region, which is a region within the tissue slab where $D  = 0$ in Equation (\ref{eq:model}). In practice, this region is never electrically activated, such that we can model it as a hole in the domain, $\Omega_{\text{non-conducting}}$. In the rest of the domain, we reduce the diffusion tensor to $D_{\text{unhealthy}} =\gamma D_{\text{healthy}}$ to account for conduction slowing $\gamma$. By default, we set $\gamma=0.8$, but we also conduct parameter estimation with $\gamma=0.1$. 

To simulate the reentry that occurs during IART, we impose a spiral initial condition on the monodomain model. Section \ref{prepacing} explains how this spiral initial condition is generated through a carefully set up prepacing simulation. Starting from this initial condition, the model is solved over a time interval $[0,T_{\text{experiment}}]$ with $T_{\text{experiment}} = 2000~\text{ms}$ for experiments with $\gamma=0.8$ and $T_{\text{experiment}} =3000~\text{ms}$ for experiments with $\gamma=0.1$  such that five periods of the spiral are simulated. Parameter inference is done based on the last two reentries, as we consider them to be the steady state of the spiral. Therefore, the model prediction is based on the simulation interval $[T_0,T_{\text{experiment}}]$ with $T_0$ the time at which the wave passes for the first time at the electrode with index 4 during the one-to-last reentry. 

Concretely, we replace $V_m(x,t)$ by $V_m(x,t,\vartheta)$ in Equations (\ref{eq:model}) and (\ref{eq:h}) to denote its dependency on this parameter value and solve them
\begin{align}
    x &\in \Omega=\Omega_{\text{atrium}} \setminus \Omega_{\text{non-conducting}} &&t\in[0,T_{\text{experiment}}],
    \end{align}
and for the initial condition 
\begin{align} \label{eq:init}
    V_m(x,0)= V_{m,0}
\end{align}
 Figure \ref{fig:pers-model} shows an example of this type of domain, together with the initial condition.

\subsection{Geometric parameter inference}\label{personalized-IART}
In this work, we calibrate the model of cardiac EP during IART from Section \ref{IART-model} to data by inferring a geometric parameter that describes the non-conducting region $\Omega_{\text{non-conducting}}$. We approximate this region as an ellipse characterized by its two radii, $a$ and $b$, and its inclination angle, $\varphi$ (see Figure~\ref{fig:pers-model}). Inference of this region becomes estimation of the geometrical parameter \begin{equation}\label{eq:parameter} 
    \vartheta=[a,b,\varphi],
\end{equation}
and we obtain a parameterized description of the non-conducting region
\begin{align}
    \Omega_{\text{non-conducting}} (\vartheta) = \{(x,y)\in \mathbb{R}^{2}:&
     \frac{((x-x_0)\cos(\phi)+(y-y_0)\sin(\phi))^{2}}{a^{2}}\\ \nonumber
     \qquad& +\frac{((x-x_0)\sin(\phi)+(y-y_0)\cos(\phi))^{2}}{b^{2}}<1\}
\end{align}
 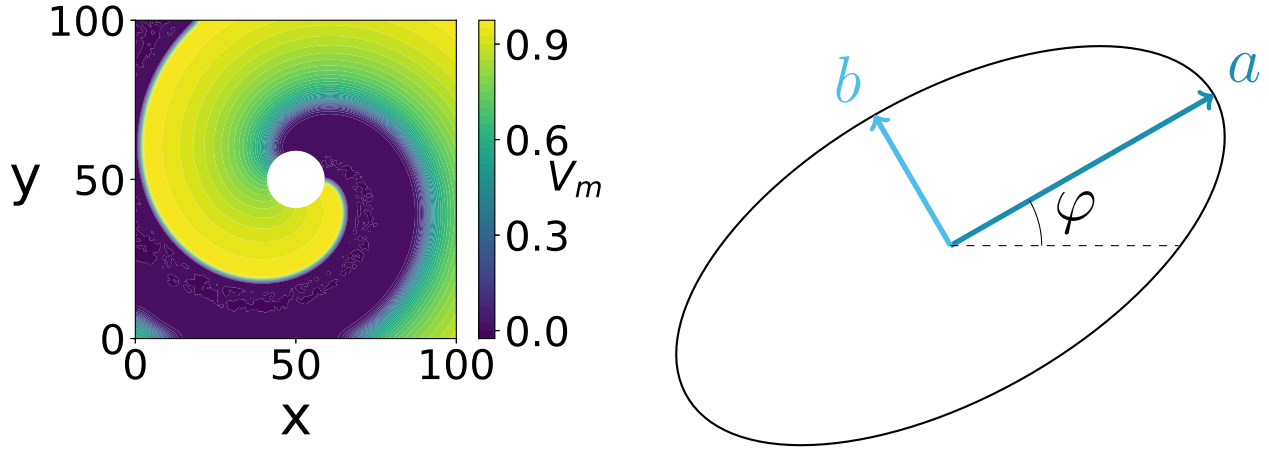
\begin{figure}[ht]
     \centering
\begin{minipage}[t]{0.49\linewidth}
	\centering
	\includegraphics[width=\linewidth]{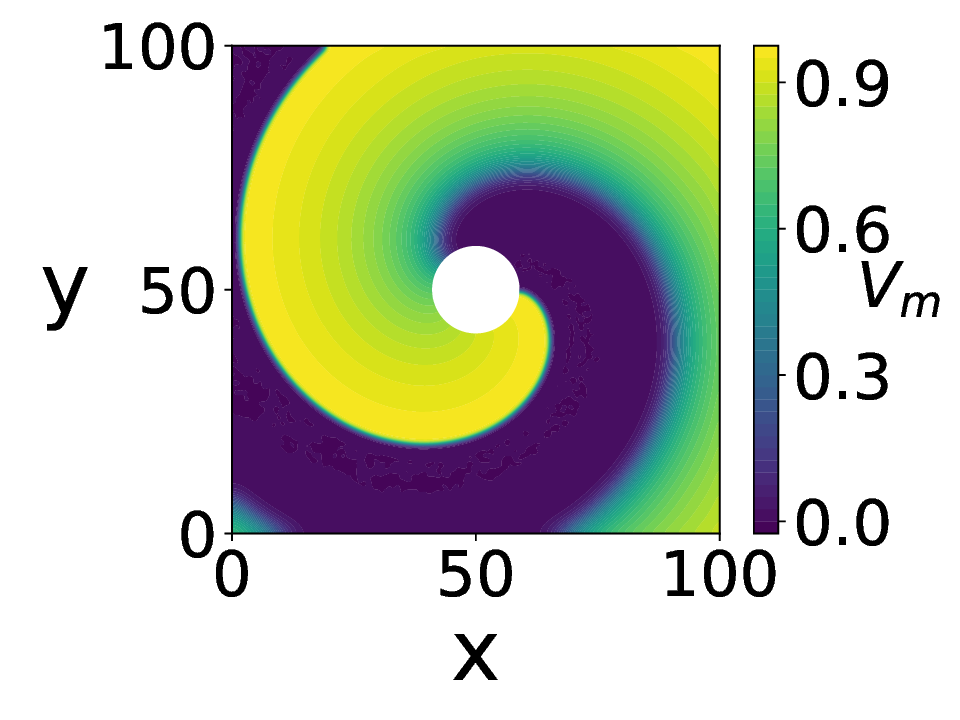}

\end{minipage}
\hfill
\begin{minipage}[t]{0.49\linewidth}
	\centering
	\definecolor{KULeuvenLichtblauw}{cmyk}{.59,.07,.00,.00}
\definecolor{KULeuvenDeptCS}{cmyk}{.65,.00,.00,.35}

\begin{tikzpicture}
  \def\phi{30} 

  \begin{scope}[rotate=\phi]
    \draw[thick] (0,0) ellipse (4 and 2);
  \draw[->, KULeuvenDeptCS, line width=2pt] (0,0) -- ({4},{0}) node[above right] {\Huge $a$};

  \draw[->, KULeuvenLichtblauw,line width=2pt] (0,0) -- (0,2) node[above left] {\Huge $b$};
  \end{scope}

  \draw[dashed] (0,0) -- (3,0);

  \draw (1.2,0) arc[start angle=0,end angle=\phi,radius=1.2];
  \node at ({1.7*cos(\phi/2)},{1.7*sin(\phi/2)}) {\Huge $\varphi$};

\end{tikzpicture}
\end{minipage}
\caption{Left: the IART behavior occurs inside a 2D tissue slab with a non-conducting elliptical region, with $D=0$ in equation \ref{eq:model} in this region, and a spiral wave as the initial condition. Right: the geometrical parameter that describes the non-conducing region. The plot on the left takes $\vartheta=[9.0~\text{mm},9.0~\text{mm},0.0~\text{rad}]$}
     \label{fig:pers-model}
 \end{figure}

The methods proposed in this paper can be directly applied to other parametric representations of the non-conducting region of low or moderate dimension. A natural extension would be to include the location of the region’s center as an additional parameter. However, this study fixes this location to avoid practical complications that arise when trying out different values for the parameter $\vartheta$ during parameter estimation causes measurement points to coincide with the non-conducting region (see Section \ref{data} and \ref{Bayes}). Our low-dimensional description of this non-conducing region fits into our strategy to increase the difficulty of the inference problem stepwise. We start from an initial toy example that will be used as a base case to extend our methods to more realistic settings in follow-up work aimed at mitigating the growing computational burden associated with inference of high-dimensional parameters.

 \subsection{Numerical solution}\label{numerical-solution}
To numerically solve Equation (\ref{eq:model})--(\ref{eq:init}), we need to specify a spatial and time resolution, for which we introduce the notation $\Delta = \{\Delta x,\Delta t\}$. In our experiments, we set the resolution to $\Delta  =\{0.25~\text{mm}, 0.5~\text{ms}\}$ by default. However, for some experiments, the spatial discretization will be coarsened to $\Delta x =0.5~\text{mm}$ to investigate the impact of the discretization error. In space, we discretize the 2D tissue slab and generate a mesh using GMSH \cite{geuzaine2009} with a targeted element size of $\Delta x$. The mesh models the non-conducting region specified by $\vartheta$ as a hole in the domain. 

To solve the monodomain over this mesh, we use the software package fenicsx-beat\footnote{https://github.com/finsberg/fenicsx-beat} that implements the splitting solver for the monodomain model \cite{sundnes2006}. The cell model is downloaded from the CellML repository \cite{cuellar2003} and integrated into fenicsx-beat. 

The splitting solver method decouples the diffusion and the reaction term into a spatially coupled system of PDEs and a spatially decoupled system of ordinary differential equations (ODEs), respectively. It solves them sequentially \cite{sundnes2006} by taking the following steps for every time step $\Delta t$: 
\begin{enumerate}
    \item The reaction term is solved over $\frac{\Delta t}{2}$.
    \item The result from step 1 is used as the initial condition to solve the diffusion term over a full time step $\Delta t$.
    \item The result from step 2 is then used to solve the reaction term again over $\frac{\Delta t}{2}$.
\end{enumerate}
For a more detailed description of this scheme, we refer the reader to \cite{sundnes2006}.

The contribution of the reaction term involves solving a system of ODEs (one for each mesh point), for which we employ the first-order Rush-Larsen scheme \cite{sundnes2006}. The diffusion term requires a PDE solve with the Finite Element Method (FEM), for which the FEM provided by the FEniCS project is used \cite{alnaes2013,baratta2023, scroggs2022,scroggs2022a}. We use Lagrange elements with first-order basis functions. For a spatial resolution given by $\Delta x=0.25~\text{mm}$, this results in approximately ${220\,000}$ degrees of freedom.
Our result is a numerical approximation for $V_m(x,t,\vartheta)$
\begin{equation} \label{eq:discretized monodomain}
    [\hat{V}_m(x,t_{m,\Delta}(\vartheta),\hat{V}_{m,\Delta}(x,t_1,\vartheta),  \dots,\hat{V}_{m,\Delta}(x,t_{N_{sim}},\vartheta)]
\end{equation}
with $t_i=T_0+i\Delta t$ and $N_{\text{frame}}=\frac{T_{\text{experiment}}-T_0}{\Delta t}$. 

\subsection{Prepacing simulation} \label{prepacing}
To derive the spiral initial condition, we set up a prepacing simulation of the monodomain model that performs an S1-S2 protocol. Such a protocol stimulates two regions in the heart sequentially. When carefully designed, part of the tissue is still recovering from the first stimulus (the S1-stimulus) at the moment of the second stimulus (the S2-stimulus), such that the resulting wave can only propagate in one direction and generates a spiral wave.  The motivation behind the prepacing simulation is reducing computational time, as we now only simulate the pacing that creates the reentrant tachycardia once, instead of during each simulation. 

The prepacing simulation solves the monodomain equation over a two-dimensional reference mesh with a circular non-conducting region centered at $(5\;\text{cm},5\;\text{cm}) $ and with a radius of $9$~mm, which is equivalent $\vartheta = [9.0~\text{mm},9.0~\text{mm},0.0~\text{rad}]$. The conduction slowing is set to $\gamma = 0.8$. We run this simulation over a long time interval of $[0,T_{\text{prepacing}}]$ with $T_{\text{prepacing}} = 100\,000$~ms to ensure that the resulting spiral is in steady state. This prepacing is performed numerically following the setup from Section \ref{numerical-solution}. As will become clear in Section \ref{Bayes}, this chosen non-conducting region is an average over all its configurations that are considered during inference, which ensures that the imposed initial condition is close to the steady-state reentry for each considered configuration and the required duration is reduced for each simulation.

\section{Modeling of synthetic measurement data} \label{data}
After defining the forward model, this section presents the synthetic data that serve as input for the Bayesian estimation of the non-conducting region. First, Section \ref{pentaray} describes the mapping catheter that measures the EGM data during clinical procedures. Second, we outline how we synthetically generate data that mimic the measurements of this catheter in Section \ref{setup}.

\subsection{Mapping catheter} \label{pentaray}
\begin{figure} [ht]
    \centering
    \includegraphics[width=0.45\textwidth]{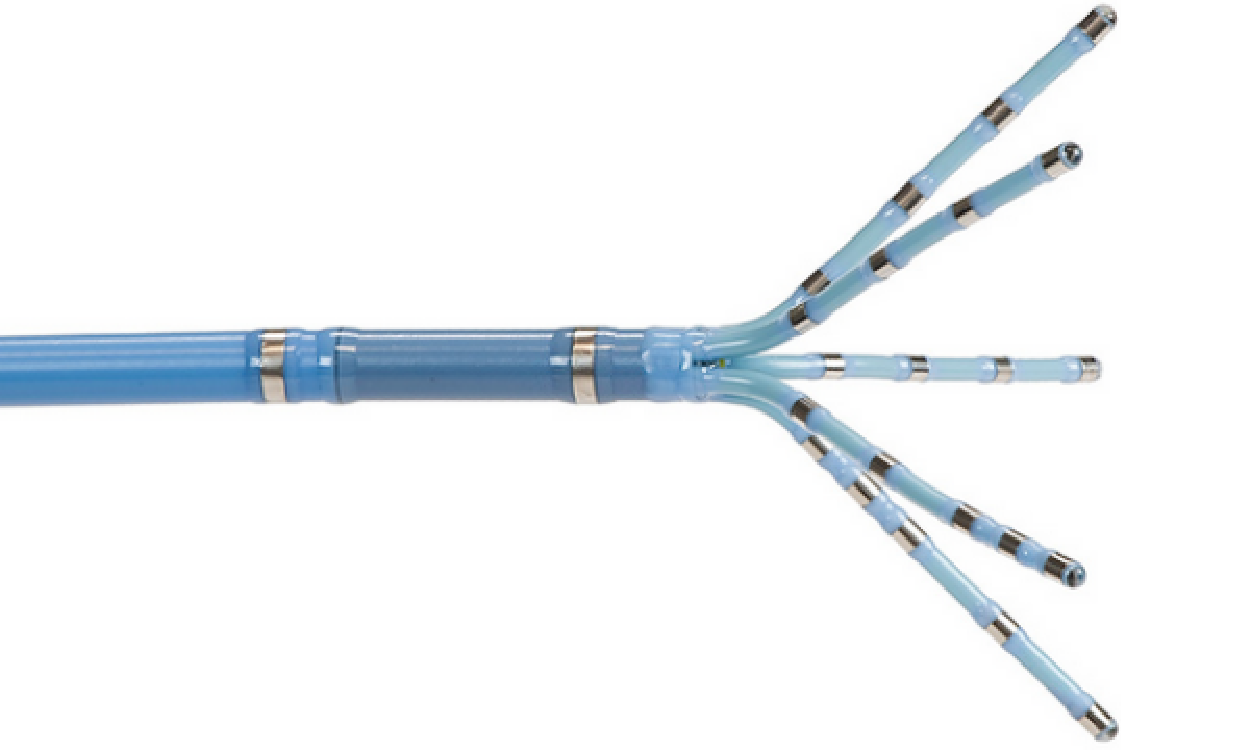}
    \includegraphics[width=0.39\textwidth]{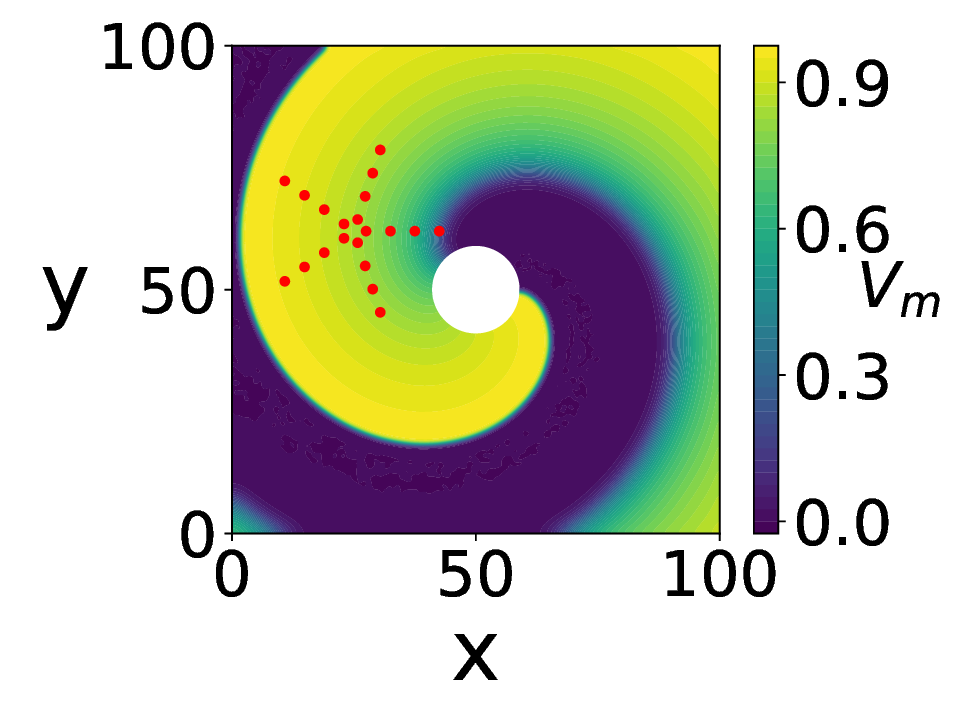}
    \todo[inline]{Change figure b to fenics}
    \caption{Left: the $\text{PENTARAY}^\text{®}$ NAV ECO High Density Mapping Catheter (Biosense Webster, Diegem, Belgium) consists of 5 arms with 4 electrodes (grey) on each. Right: the location of the electrodes is projected on the 2D plane of the tissue. Endocardial measurements are made at the blue locations for all experiments, with additional measurements at the red locations for the experiment with two catheter locations.}
    \label{fig:data}
\end{figure}

We perform inference using synthetic EGM data, generated to mimic the measurements obtained with a $\text{PENTARAY}^\text{®}$ NAV ECO High Density Mapping Catheter (Biosense Webster, Diegem, Belgium). This catheter is used for invasive mapping and consists of five flexible arms, each equipped with four electrodes spaced 4 mm apart (see also Figure \ref{fig:data}). In this work, the catheter is assumed to be at a fixed location at a distance from the non-conducting region. All electrodes are in direct contact with the endocardial surface and record the local electrical field. The exact electrode locations are visually illustrated as black dots in Figure \ref{fig:data} and are detailed in the supplementary material. We only consider electrode locations at a distance from the non-conducting region to guarantee local electrical activity at each electrode. This will be important for our likelihood definition based on local electrical activity in Section \ref{Bayes}.

\subsection{Generation of synthetic data} \label{setup}
We generate synthetic data to resemble the EGM recordings during IART. First, we solve the IART model described in Section \ref{IART-model} with an elliptical non-conducting region defined by the true geometrical parameter $\vartheta_{\text{true}} = [10.0\text{ mm} ,4.0\text{ mm},0.0\text{  rad}]$ and solve the model numerically as outlined in \ref{numerical-solution} to obtain $\hat{V}_{m,\Delta}(x,t_i,\vartheta_{\text{true}})$ Once more, the simulation is performed over $T_{\text{experiment}} = 2000~\text{ms}$ such that five spirals are simulated.  

Then, we obtain the EGM-data by mapping the transmembrane potential to the potential measurements that would be recorded by the catheter, using the pseudo-EGM formulation \cite{bishop2011}. For a measurement at electrode location $r_e$ at time $\tau_k$, we can calculate the potential as
\begin{align}\label{eq:bidomain}
  \phi_e(r_e,\tau_k,\vartheta_{\text{true}}) &= \frac{1}{4\pi\sigma_b}\int_\Omega \frac{\nabla \cdot \sigma_i \nabla \tilde{V}_{m,\Delta}(r,\tau_k,\vartheta_{\text{true}})}{||r_e-r||}d\Omega\\
	&= \frac{1}{4\pi \sigma_b}\int_\Omega \sigma_i \nabla \tilde{V}_{m,\Delta}(r,\tau_k,\vartheta_{\text{true}}) \cdot \frac{r_e-r}{||r_e-r||^3}d\Omega,
    \end{align}
    where $\sigma_b = 1.0$~mS/mm and $\sigma_i = 0.174$~mS/mm, accounting for the blood and the intracellular conductivity, respectively. The spatial variable $r$ is three-dimensional. It differs from the two-dimensional spatial variable $x$ used in Equations (\ref{eq:model})--(\ref{eq:init}), which is only defined in the tissue. Inside the tissue, the relation between both coordinates is given by $r=[x,0]$. Moreover, we now have a transmembrane potential function in terms of the new space coordinates $\tilde{V}_{m,\Delta}(r,t,\vartheta_{\text{true}})=\hat{V}_{m,\Delta}(r_{0:2},t,\vartheta_{\text{true}})$ if $r\in \Omega$ and $\tilde{V}_{m,\Delta}(r,t,\vartheta_{\text{true}})=0$ otherwise. The first equality is the original pseudo-EGM formulation \cite{coveney2020}. Using the divergence theorem, we can rewrite it so that only first-order derivatives occur. This allows us to use the first-order basis functions in the finite element discretization and to approximate the integral with a quadrature formula. 
    
We assume that a potential measurement is made every 4~ms. Concretely, this gives, as data $y$, the time traces
\begin{equation}\label{eq:y}
    y = \begin{bmatrix}
      \phi_e(r_{1},\tau_0)&\phi_e(r_{1},\tau_1)&\phi_e(r_{1},\tau_2)&\dots& \phi_e(r_{1},\tau_{N_{\text{frame}}})\\
      \phi_e(r_{2},\tau_0)&\phi_e(r_{2},\tau_1)&\phi_e(r_{2},\tau_2)&\dots& \phi_e(r_{2},\tau_{N_{\text{frame}}})\\
     &\dots&&\ddots& \dots\\
      \phi_e(r_{20},\tau_0)&\phi_e(r_{20},\tau_1)&\phi_e(r_{20},\tau_2)&\dots& \phi_e(r_{20},\tau_{N_{\text{frame}}})\\
\end{bmatrix} + \eta,
\end{equation}
where $r_{j}$ is the coordinate of electrode $j$ of the catheter, $\tau_k = T_0 + 4k~\text{ms}$ and the number of frames $N_{\text{frame}}=\frac{T_{\text{experiment}}-T_0}{4~\text{ms}}$. We notice that we only need to map the $\hat{V}_{m,\Delta}(x,t_i,\vartheta_{\text{true}})$ from Equation (\ref{eq:discretized monodomain}) to the EGM if there is a $k$ such that $t_i=\tau_k$, which happens only once per eight time steps. To account for measurement error, we add the matrix $\varepsilon$ with $\varepsilon_{jk} \sim \mathcal{N}(0,\sigma^2)$, yielding independent Gaussian measurement noise for each measurement with variance $\sigma^{2}=10^{-6}$.

\section{Bayesian inversion problem}\label{Bayes}
We define the estimation of the geometrical parameter describing the non-conducting region from the synthetic electrogram data as a Bayesian inversion problem. This probabilistic formulation not only provides an estimate but also a quantitative measure of the uncertainty arising from measurement error. First, we recall the general definition of Bayesian inference from Equation (\ref{eq:bayesianinference}) and specify all its ingredients in Section \ref{bayes-inference}. Then, we bring the likelihood function into its final form in two steps. We propose a compressed likelihood based on characterizing quantities of the data in Section \ref{summary-statistic}, which improves the sampling efficiency and the accuracy of the noise model. We establish the impact of the discretization error on the continuity of the forward model and account for this error by inflating the variance of the discretized compressed likelihood in Section \ref{disc-ll}. Finally, Sections \ref{summary statistics calculation} and \ref{Covariance} explain the practical details of the calculation of the characterizing quantities and of the covariance of the discretized compressed likelihood function, respectively. 
\subsection{Bayesian framework} \label{bayes-inference}

We first specify the elements of Equation (\ref{eq:misfit}). Bayesian inference models the parameter $\vartheta$, the data $y$ and the measurement error $\varepsilon$ as random variables. In our setting, the parameter $\vartheta$ is the geometrical parameter given by Equation (\ref{eq:parameter}) and the data $y$ are the time traces presented in Equation (\ref{eq:y}). The map between $\vartheta$ and $y$ consists of the solution of our forward problem $\mathcal{F}(\vartheta)$, given by equations (\ref{eq:model})--(\ref{eq:init}), followed by application of the observation map $h(.)$ to the solution to obtain a matrix $H$ of endocardial potentials using the formalism from Equation (\ref{eq:bidomain}). For $h_{jk}$ the element of $H$ on row $j$ and column $k$, then $h_{jk}=\phi_e(r_{j},\tau_k$ with $i=1,2,...N_{\text{frame}}$ and $j=1,2,\dots,20$, respectively. The value of $r_{j}$ is the coordinate of electrode $j$  and $\tau_k=T_0+4k$. 

Then, we specify the elements of Equation (\ref{eq:bayesianinference}), which gives an expression to find the posterior distribution $P(\vartheta|y)$ on the geometrical parameter conditional on the data. We define a prior distribution $\pi(\vartheta)$ that captures all prior knowledge about the geometrical parameter. In our experiments, we assume no prior knowledge by defining $p(\vartheta)$ as a wide uniform distribution $a,b\sim\mathcal{U}(2~\text{mm},16~\text{mm})$ and $\varphi\sim\mathcal{U}(-\pi/2,\pi/2)$. The other main ingredient is the likelihood $\mathcal{L}(y|\vartheta)$, which gives a score of how well the model can predict the data. The definition of the likelihood requires a more elaborate description and is postponed to Sections \ref{summary-statistic} and \ref{disc-ll}.

\subsection{Compressed likelihood based on characterizing quantities} \label{summary-statistic} In the likelihood, we quantify the probability that the discrepancy between the model output and the data is attributable to measurement noise. Formally, it is defined as \begin{equation}
    \mathcal{L}(y|\vartheta) = P(y|\vartheta) = f_{\varepsilon}(y-h(\mathcal{F}(\vartheta))).
\end{equation}
with $f_{\varepsilon}(.)$ the probability distribution function of the measurement error $\varepsilon$. The likelihood assigns a score to each parameter $\vartheta$, indicating how likely the corresponding non-conducting region would have produced the observed data $y$ given the model for the measurement error. As shown in Equation (\ref{eq:bayesianinference}), the likelihood function uses these scores to update the prior distribution into a posterior distribution that incorporates the observed data $y$.

In this work, we reduce the dimension of the data from 20 time traces to a lower-dimensional vector $s_y$ of characterizing quantities, such that Equation (\ref{eq:misfit}) becomes
\begin{equation}\label{eq:misfit-ss}
    s_y = s(\mathcal{F}(\vartheta)) + \varepsilon'.
\end{equation}
The observation operator $h(.)$ has been replaced by a new map $s(.)$ that calculates the characterizing quantities from the model solution $V_m(x,t,\vartheta)$. Moreover, we now define $\varepsilon'\sim\mathcal{N}(0,\Sigma_{\varepsilon'})$. as the discrepancy between these characterizing quantities for the model prediction and the data. We define the compressed likelihood as
\begin{equation}
    \mathcal{L}(s_y|\vartheta) = P(s_y|s_\vartheta) = f_{\varepsilon'}( s_y - s(\mathcal{F}(\vartheta))).
\end{equation}
with $f_{\varepsilon'}(.)$ now the probability distribution function of the discrepancy $\varepsilon'$ in the characterizing quantities.

The dimensionality reduction paves the way for accurate inference and uncertainty quantification. After all, a good likelihood requires an accurate model for the measurement noise, and deriving such a model from clinical data is more achievable for these characterizing quantities. The primary trade-off of this approach is a potential loss of information. However, appropriate selection of the characterizing quantities should limit this loss.

We define two types of characterizing quantities: the period of the spiral around the non-conducting region and the local activation time (relLAT) at each of the electrodes, relative to the average activation time over all electrodes. An interesting benefit of these quantities is their tce of the initial condition and invariance to temporal shifts.  Since we compare the data of the steady-state spiral with the model output during steady-state behavior, this allows direct comparison of the model output to the observed data, without requiring time trace alignment. Moreover, it paves the way for using experimental data where the initial condition is often unknown. 

This results in the following expression for the compressed likelihood function
\begin{equation} \label{eq:ll}
    \mathcal{L}(s_y|\vartheta)=  \frac{1}{(2\pi)^{d_s/2}\det(\Sigma_{\varepsilon'})^{1/2}} \exp\left( -\frac{1}{2}  (s_y - s(\mathcal{F}(\vartheta)))^T\Sigma_{\varepsilon'}^{-1}(s_y - s(\mathcal{F}(\vartheta)))\right),
\end{equation}
where $d_s$ is the dimension of the vector of characterizing quantities. The calculation of $s_y$ from the data $y$ and the operator $s(.)$ are detailed in Section \ref{summary statistics calculation}. The derivation of the covariance  matrix $\Sigma_{\varepsilon'}$ follows in Section \ref{Covariance}

While we adopt a simple Gaussian noise model and a limited set of quantities in this study, the proposed methodology is compatible with more accurate noise models and richer sets of characterizing quantities. These inputs would influence the posterior distribution, but not the underlying inference framework or the methods to find the posterior distribution. As this paper focuses on methodological developments rather than modeling clinical uncertainties, the incorporation of more detailed clinical modeling is left for future research. 

\subsection{Discretized likelihood evaluation} \label{disc-ll}
The compressed likelihood function involves an evaluation of the forward model $\mathcal{F}(\vartheta)$ (see Equation \ref{eq:ll}). As this requires solving a PDE, we will have to settle for a numerical approximation. Following the approach laid out in Section \ref{numerical-solution}, we approximate $\mathcal{F}(\vartheta)=V_m(x,t,\vartheta)$ as 
 \begin{equation}
     \hat{\mathcal{F}}^{(\text{IM})}_{\Delta}(\vartheta)=[\hat{V}_{m,\Delta}(x,t_0,\vartheta),\hat{V}_{m,\Delta}(x,t_1,\vartheta),  \dots,\hat{V}_{m,\Delta}(x,t_{N_{sim}},\vartheta)],
 \end{equation}
 with $\Delta = \{\Delta x, \Delta t\}$,  $t_i=T_0+i\Delta t$ and $N_{sim}=\frac{T_{\text{experiment}}-T_0}{\Delta t}$. We introduce the superscript (IM), which is used as an abbreviation for Independent Meshing. For consistency over all plots in this paper, we introduce this notation while we postpone its explanation to later in this section.

\begin{figure}[ht]
     \centering
    \includegraphics[width=0.45\textwidth]{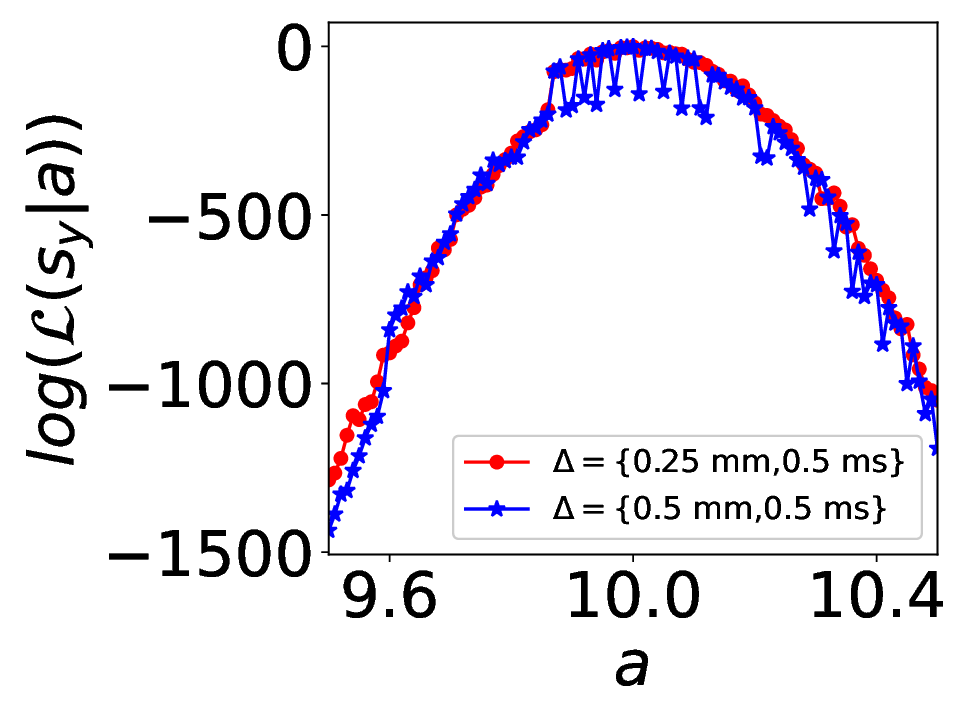} 
    \includegraphics[width=0.45\textwidth]{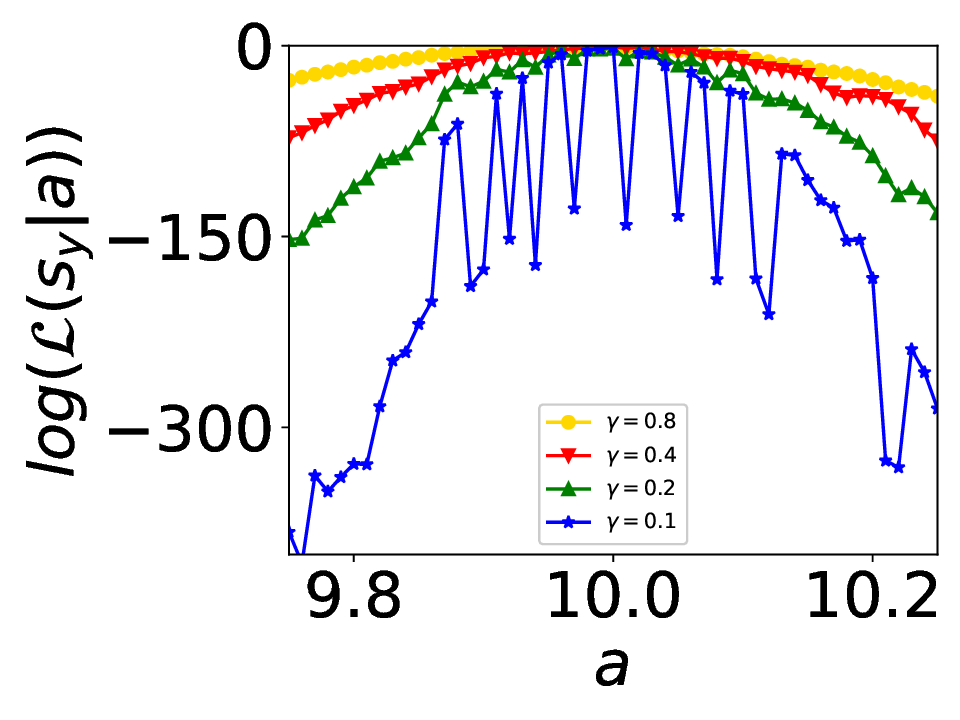} 
     \caption{Plots of the compressed likelihood as a function of the long radius $a$ for an inference experiment that only estimates $a$, with ground truth $a_{true}=10.0$~mm. We set the conduction slowing to $\gamma=0.1$ and the likelihood variance to $\Sigma_{\varepsilon'}=\diag(0.1,1.0,1.0,\dots,1.0)$. Left: the compressed likelihood using different resolutions $\Delta x$. Right: the discretized compressed likelihood for different $\gamma$ with $\Delta x=0.5$~mm fixed.}
     \label{fig: discontinuous ll}
 \end{figure}

The discretization error causes an extra discrepancy $d =s(\mathcal{F}_{\Delta}^{(\text{IM})}(\vartheta))-s(\hat{\mathcal{F}}(\vartheta))$ into Equation (\ref{eq:misfit-ss})
 \begin{equation}
     s_y = s(\hat{\mathcal{F}}^{(\text{IM})}_{\Delta}(\vartheta))+\varepsilon' + d.
 \end{equation}
 In Bayesian inference studies, it is common to plug the discretized forward model $\hat{\mathcal{F}}^{(\text{IM})}_{\Delta}(\vartheta)$ directly into the likelihood without further consideration of $d$. However, for accuracy, this must be preceded by a formal check of the implied assumption $\varepsilon'\gg d$. 
 Otherwise, we risk underestimating the discrepancy between model and data, resulting in a posterior that is overconfident in incorrect estimates with significant associated risks \cite{poot2025}. In our inference problem with a geometrical parameter, we observe an even more pronounced effect: the posterior has artificial discontinuities solely due to the underestimation of the discretization error. Figure \ref{fig: discontinuous ll} illustrates these discontinuities with a plot of a numerical approximation to the compressed likelihood from Equation (\ref{eq:ll}) as a function of the discretization $\Delta x$. Logically, discontinuities in the compressed likelihood are transferred into the posterior, where they can act as visual warnings for underestimation of the discretization error. While this could be seen as a complication, the advantage is that these discontinuities act as an a posteriori visual warning for discretization error underestimation. 

The remainder of this section now addresses the three following matters. First, we explain why the likelihood is inherently discontinuous. Then, we discuss why the assumption $\varepsilon' \gg d$ is more easily violated in our setting, which in turn causes these discontinuities to be large, such that they have an impact on the posterior shape. Finally, we propose a methodology to proceed when this assumption is violated.

\begin{figure}[ht]
     \centering
     \input{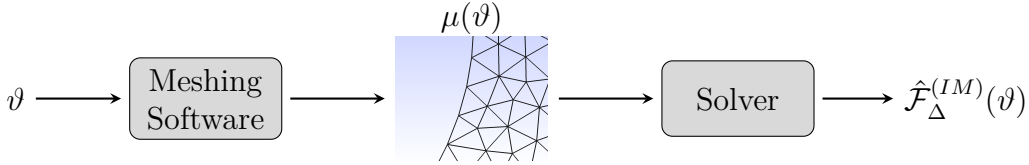}
     \caption{In shape estimation problems, the numerical solution of the forward model for the shape $\vartheta$ consists of a discontinuous meshing step to generate a mesh $\mu$, that models the shape $\vartheta$, and the FEM step to solve the discretized model over this mesh.}
     \label{fig:shape estimation}
 \end{figure}

To explain the discontinuities, Figure \ref{fig:shape estimation} shows the steps to solve the forward model $\hat{\mathcal{F}}^{(\text{IM})}_{\Delta}(\vartheta)$ for a value of $\vartheta$.  First, we generate a mesh $\mu(\vartheta)$ that models the region described by $\vartheta$. Then, we solve the discretized PDE numerically over this mesh to obtain $\hat{\mathcal{F}}^{(\text{IM})}_{\Delta}(\vartheta)$. Even for small perturbations to $\vartheta$, the meshing step is done completely from scratch. As a result, the meshing step is discontinuous in function of $\vartheta$ and makes the overall discretized forward model solve $\hat{\mathcal{F}}^{(\text{IM})}_{\Delta}(\vartheta)$ discontinuous as well, even if the original forward model $\mathcal{F(\vartheta)}$ is not. In turn, this causes the compressed likelihood and the posterior to be discontinuous as well. In the remainder of this paper, we call this meshing approach Independent Meshing to contrast it with an alternative approach proposed in Section \ref{MCMC-ar}.

In our setting, the effect of the discretization error is amplified by the observable map $s(.)$. We illustrate this using Figure \ref{fig: discontinuous ll}, which shows that, besides the resolution $\Delta x$,  the discretization error depends on the value of the conduction slowing $\gamma$ through its impact on the diffusion constant $D$ (i.e., $D= \gamma D_{\text{healthy}}$) in Equation (\ref{eq:model}) as well. Slower waves increase all the characterizing quantities in our likelihood. The same relative error due to the discretization will be larger in absolute terms and, as a result, the relative importance of $d$ in comparison to $\varepsilon'$ increases. While $\varepsilon'\gg d$ still holds for $\gamma=0.8$, it is no longer the case for $\gamma=0.2$ and $\gamma=0.1$. This observation reinforces our claim that the size of the discretization error should be formally checked. After all, the discretization of mathematical EP models is typically chosen based on benchmarking studies \cite{niederer2011} or requirements for modeling spiral waves. While a resolution of $\Delta x=0.5$~mm might be fine enough to obtain a good numerical approximation of the overall solution (as given by Equation (\ref{eq:discretized monodomain})) and even for the characterizing quantities for some values of $\gamma$, we might have incorrectly assumed $\varepsilon'\gg d$ for other model specifications.  
 
When $\varepsilon'\gg d$ does not hold for an experiment, we model the discretization error in the discretized compressed likelihood rather than enforcing it by refining the resolution. After all, sufficiently increasing the resolution would become completely infeasible for small values of $\gamma$. Recently, a review \cite{poot2025} compared different approaches to do so, of which only the statistical finite element method (statFEM) \cite{girolami2021} and the random mesh finite element method (RM-FEM) \cite{abdulle2021} can be applied to our inference problem. We employ a method derived from statFEM and model the discretization error as a random variable that we assume to be Gaussian
 \begin{equation}
     d \sim \mathcal{N}(0,\Sigma_d).
 \end{equation}
This variance can be estimated from the data during inference as a hyperparameter, but this requires highly informative data. Therefore, we opt to estimate $\Sigma_d$ before inference, following the procedure explained in Section \ref{Covariance}. In this work, we use the same numerical precision for our model as for the generation of the synthetic data. This ensures that the discretization-induced bias is effectively zero, allowing us to focus solely on addressing the problem-specific challenge posed by discontinuities in the likelihood via likelihood variance inflation. However, statFEM also provides an additional parameter for modeling discretization bias. Consequently, our framework could be extended to incorporate a more realistic setting where this bias is non-zero.

The discretized compressed likelihood becomes
 \begin{equation}\label{eq:ll-final}
 \begin{split}
\hat{\mathcal{L}}^{(\text{IM})}_{\Delta }(s_y|\vartheta)=  &\frac{\exp\left( -\frac{1}{2} (s_y - s(\hat{\mathcal{F}}^{(\text{IM})}_{\Delta }(\vartheta)))^T\Sigma^{-1}(s_y - s(\hat{\mathcal{F}}^{(\text{IM})}_{\Delta }))\right)}{(2\pi)^{d_s/2}det(\Sigma)^{1/2}}\\ 
\end{split}
\end{equation}
with $\Sigma=\Sigma_{\varepsilon'}+\Sigma_d$. While this improves the accuracy of the inference by modeling the discretization error, the discontinuities in the likelihood and posterior remain, albeit in reduced size, and their influence on the posterior shape is damped. 

\subsection{Calculation of characterizing quantities}\label{summary statistics calculation} In this section, we explain how we can calculate the characterizing quantities of the data $y$ and the model solution $\mathcal{F}(\vartheta)$. To derive the characterizing quantities, we first need to calculate the local activation time (LAT) at each electrode. We introduce the notation $\widehat{LAT1}^{\alpha}_j$ and $\widehat{LAT2}^{\alpha}_j$ to denote the arrival time at electrode j during the first and second reentry, respectively. The superscript $\alpha$ is either $\alpha=y$ or $\alpha=\mathcal{F}(\vartheta)$, to annotate whether a particular LAT is derived from the data or the model prediction, respectively. In both situations, the meaning of LAT is equivalent, but the procedure to calculate it differs strongly. This procedure is not described in the core of the paper for brevity, but is included in the Supplementary Material for reproducibility purposes.

From the LATs, we compute the relLAT at an electrode as the difference between the LAT at that electrode and the average across all electrodes for the last full reentry within the simulation interval,  $\widehat{LAT}2^{\alpha}$. We calculate the period as the average over the LAT intervals between the last two full reentries at all electrodes. We thus have:
\begin{align} \label{eq:stats-calc-1}
    \widehat{LAT1}^{\alpha}&=\frac{1}{20} \sum^{20}_{j=1} LAT1^{\alpha}_{j}\\
    \widehat{LAT2}^{\alpha}&=\frac{1}{20} \sum^{20}_{j=1} LAT2^{\alpha}_{j}\\
    relLAT_{j} &= LAT2^{\alpha}_{j} - \widehat{LAT2}^{\alpha}\\
    \label{eq:stats-calc-2}
    period &=  \hat{LAT2}^{\alpha}- \widehat{LAT1}^{\alpha}.
\end{align}

\subsection{Covariance matrices} \label{Covariance}
In this section, we specify the covariance matrix $\Sigma$, which contains one contribution of the measurement noise, $\Sigma_{\varepsilon'}$, and one contribution of the discretization error, $\Sigma_d$, to fully define the likelihood in Equation (\ref{eq:ll-final}). For the former, the derivation of a noise model based on clinical data is ongoing work. In this work, we define a simplified matrix $\Sigma_{\varepsilon'}$ based on the assumption that the error on each LAT measurement has a standard deviation of $1.0$~ms. This error includes both the measurement error and the errors made to derive the LAT from the electrogram measurements. Under the assumption that the measurement error at different electrodes is independent, application of Equations (\ref{eq:stats-calc-1})--(\ref{eq:stats-calc-2}) to normal random variables results in measurement noise with variance $0.1$ for the period and $1.0$ for each relLAT and $\Sigma_{\varepsilon'}=\diag(0.1,1.0,\dots,1.0)$. 

We estimate $\Sigma_d$ on a case-by-case basis. The variance of each characterizing quantity due to the discretization error is calculated by solving the model for $\vartheta_{\text{true}}$ on 51 different mesh configurations. To this end, we first create a mesh for $\vartheta=9.95 + 0.02s$ with $s=0,\dots,51$, and then, we relocate its nodes such that this mesh fits $\vartheta_{\text{true}}$. We calculate each characterizing quantity's sample variance over all model solutions and multiply these variances by 1.3, a heuristic value, to safely extrapolate this variance to other values of $\vartheta$, in the spirit that it is better to overestimate this variance than to underestimate it. We round the resulting values to one number after the decimal point, because our heuristic estimation of these digits would not be meaningful. The resulting values are taken as diagonal elements of $\Sigma_d$.

\section{Markov chain Monte Carlo sampling}\label{MCMC}
We perform Markov chain Monte Carlo (MCMC) sampling to estimate the posterior distribution defined in Section \ref{Bayes}. In Section \ref{Bayes}, we specified the prior and the compressed likelihood that appear in Equation (\ref{eq:bayesianinference}), but the evidence remains an unknown normalization constant ensuring that the posterior is a valid probability distribution. Computing this term requires integrating the numerator over the entire parameter space, which is computationally intractable due to the complexity of the forward model embedded in the compressed likelihood. Instead, we employ MCMC to approximate the posterior distribution.

To perform MCMC sampling of the posterior distribution, we develop an adapted Metropolis-Hastings (adapted MH) algorithm to enhance the efficiency and robustness of the sampling in the posterior covariance anisotropy and discretization errors. First, Section \ref{MCMC-MH} explains the random-walk Metropolis-Hastings (RWMH) algorithm, the classical algorithm for MCMC sampling, which was taken as a starting point for our algorithm. Next, we change both the proposal and the accept-reject step of this algorithm for improved sampling. In Section \ref{MCMC-proposal}, we propose an adaptive proposal step, first proposed in the Adaptive Metropolis algorithm, for efficient exploration of the correlated parameter space. In Section \ref{MCMC-ar}, we modify the accept-reject step to account for the interaction between the meshing and the model output, which was shown in Figure \ref{fig:shape estimation}. We use node relocation to generate a mesh for the proposed sample and achieve local continuity for the discretized compressed likelihood. Combining this proposal step and accept-reject step results in our adapted MH algorithm.

\subsection{RWMH algorithm} \label{MCMC-MH}
The classical RWMH algorithm performs MCMC sampling of a posterior distribution that is known up to a constant. A Markov chain is created with the posterior as its invariant distribution. The ergodicity of the chain ensures that simulating the chain samples the invariant distribution, thus the posterior distribution. While RWMH methods are well-established, this section proceeds with a qualitative description of their setup to highlight the modifications we propose to the standard form of the algorithm.

RWMH simulates a Markov chain to generate samples of the posterior.  After $l$ iterations, we have the samples
\begin{equation}
    \{\vartheta^{(0)},\vartheta^{(1)}, \dots, \vartheta^{(l)}\}.
\end{equation}
During iteration $l+1$, we first propose a new sample based on the last sample of the posterior $\vartheta^{(l)}$
\begin{equation}\label{eq:proposal}
     \vartheta^{\ast} = \vartheta^{(l)} +u^{(l)}\;\;\;\;\;\;\;\;\;\;\;\;\;\;\;\;\;\;\;\;\;\;\text{with } u^{(l)} \sim \mathcal{N}(0,\Sigma_u).
\end{equation}
Then, this proposal is accepted with probability
\begin{equation}\label{eq:ar} 
    a(\vartheta^{(l)},\vartheta^{\ast})=\frac{P(\vartheta^{\ast}|y)}{P(\vartheta^{(l)}|y)} = \frac{\mathcal{L}(s_y|\vartheta^{\ast})\pi(\vartheta^{\ast})}{\mathcal{L}(s_y|\vartheta^{(l)})\pi(\vartheta^{(l)})}.
\end{equation}
The new sample is $\vartheta^{(l+1)}=\vartheta^{*}$, if accepted, and $\vartheta^{(l+1)}=\vartheta^{(l)}$, otherwise.

These samples can subsequently be used as an estimator of quantities of interest over the posterior. Let this quantity of interest be
\begin{equation}
    Q=\mathbb{E}[f]  = \int f(\theta) d\mu(\theta).
\end{equation}
Then, we can create the estimator
\begin{equation} \label{eq:estimator}
    \hat{Q}^{} = \frac{1}{N}\sum^{N}_{n=b} f(\theta^{(n)}),
\end{equation}
with $b$ the length of the burn-in. As we will consider warm starts, we set $b=0$.
For an appropriate definition of $f$, this offers an estimate of the mean and moments of the posterior to obtain both a parameter estimate and a quantification of the uncertainty on this estimate.

We start from the RWMH algorithm and adapt it to our case study for optimal efficiency and robustness to the discretization error. We propose an adaptation to the proposal and accept-reject step, which is detailed in Section \ref{MCMC-proposal} and \ref{MCMC-ar}, respectively.

\subsection{Proposal step}\label{MCMC-proposal}
The proposal distribution $f_u$ has an important impact on the efficiency of the sampling. Its objective is parameter space exploration, which improves if the proposal shape matches the posterior shape well. A more formal argument follows from looking at the estimator in Equation (\ref{eq:estimator}) as a Monte Carlo (MC) estimator for quantities of interest over the posterior. An MC estimator is optimal in case all samples are independent from each other. While independent samples are unachievable, the chosen proposal distribution still impacts the correlation between subsequent samples, also known as the autocorrelation, by avoiding both high rejection rates and steps that are very small compared to the posterior variance.

In RWMH, we tailor the Gaussian proposal distribution (see Equation (\ref{eq:proposal})) through adapting its covariance $\Sigma_u$ for a good exploration of the posterior distribution beforehand. Ideally, this covariance matches the covariance matrix of the posterior well, but the posterior covariance is not known in advance. Therefore, for practical reasons, we take a diagonal covariance matrix and tune its diagonal elements manually such that we achieve a good acceptance rate, ideally around $30-50\%$, which reduces the autocorrelation lengths of the Markov chain. By individually tuning the diagonal elements, we can take the anisotropy of the posterior covariance into account, as long as its principal axes coincide with individual parameters. However, a proposal with a diagonal covariance can strongly deteriorate the parameter exploration if the posterior distribution has an anisotropic covariance with principal components that do not coincide with the individual parameters. To achieve the desired acceptance rate, the manually tuned covariance will underestimate the posterior variance in certain directions and overestimate it in others. Proposal moves in directions with high variance are too small and thus have a high acceptance probability, while the opposite is true for proposals in directions with low variance. The result is a Markov chain that subsequently samples regions of the posterior, instead of sampling the entire posterior continuously. This phenomenon is called bad mixing and increases the autocorrelation of the Markov chain. As a result, longer chains are needed for convergence, requiring a larger number of samples or, equivalently,  expensive forward simulations.

The posterior correlation between the two radii $a$ and $b$ motivated us to use the proposal distribution that is adapted to such correlations \cite{haario2001}. In our algorithm, we use the adaptive Metropolis strategy, which employs the sample covariance of the previous posterior samples as the covariance of the Gaussian proposal up to a constant and thus learns the posterior structure during sampling to improve future proposals. After $l$ samples, we propose $u^{(l)}$
\begin{equation}
   u^{(l)}\sim\mathcal{N}(0,\Sigma_{u,l}),
\end{equation}
where we define the covariance matrix as
\begin{equation}\label{eq:am}
    \begin{cases}
        \Sigma_{u,l} = \Sigma_{u,0},  & \text{if $l\leq l_0$}\\
        \Sigma_{u,l} = s_d \text{cov}(\theta^{(0)},\theta^{(1)},\dots,\theta^{(l)}) +\epsilon I   & \text{if $l>l_0$ }
        
    \end{cases},
\end{equation}
with $l_0=100$, $s_d=1.152$ and $\epsilon=0.0001$. The value of $\Sigma_{u,0}$ is chosen for each experiment to obtain an acceptance rate of $30-50\%$ during the burn-in.
After a burn-in period during which the covariance structure is learned, the proposal aligns itself with the posterior, resulting in good mixing. The autocorrelation is reduced, such that the required number of samples for convergence decreases.

\subsection{Meshing strategy inside accept-reject step}\label{MCMC-ar}
Section \ref{disc-ll} explained the inherent discontinuity of the discretized forward map $\hat{\mathcal{F}}^{(\text{IM})}_{\Delta}(\vartheta)$ for a geometrical parameter describing a region and how to model the discretization error in the likelihood. This section continues this discussion and considers the impact of this error on the sampling algorithm. Concretely, we propose an adaptation to the accept-reject step for more efficient and accurate sampling in the presence of discontinuities induced by discretization errors.

 In the RWMH algorithm, the discretized forward model affects the accept-reject step through the appearance of the discretized compressed likelihood in the acceptance probability. The expression for $a(\vartheta^{(l)},\vartheta^{\ast})$ in Equation (\ref{eq:ar}) is approximated as
\begin{equation}\label{eq:ar-2} 
    \hat{a}^{(\text{IM})}_{\Delta}(\vartheta^{(l)},\vartheta^{\ast})=\frac{\hat{P}^{(\text{IM})}_{\Delta}(\vartheta^{\ast}|y)}{\hat{P}^{(\text{IM})}_{\Delta}(\vartheta^{(l)}|y)} = \frac{\hat{\mathcal{L}}^{(\text{IM})}_{\Delta}(s_y|\vartheta^{\ast})\pi(\vartheta^{\ast})}{\hat{\mathcal{L}}^{(\text{IM})}_{\Delta}(s_y|\vartheta^{(l)})\pi(\vartheta^{(l)})}.
\end{equation}
Since the discretized likelihood contains artificial discontinuities, these might have an impact on the quality of this approximation. In our algorithm, we aim to mitigate the effect of the discretization error on the accept-reject step by developing a meshing strategy tailored to Equation (\ref{eq:ar-2}).

\begin{figure}[ht]
     \centering
\input{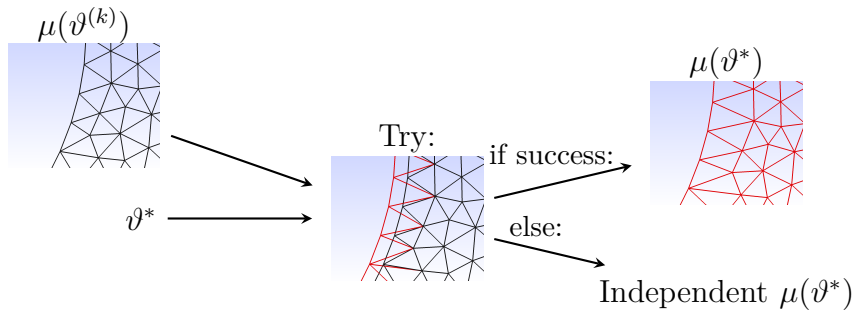}
\caption{Visual illustration of the Node Relocation strategy during iteration $k+1$: the mesh $\mu^{(k)}$ that models the last accepted non-conducting region $\vartheta^{(k)}$ is used as a starting point to generate the mesh for the proposed non-conducting region $\vartheta^{\ast}$. These nodes are relocated with the GMSH function \textit{node\_relocation} to match the non-conducting region $\vartheta^{\ast}$ with $\mu(\vartheta^{\ast})$. If the resulting mesh passes all quality checks, it is kept; otherwise, a new independent mesh is generated.}
     \label{fig:Relocating nodes}
 \end{figure}

We propose the Node Relocation strategy in Figure \ref{fig:Relocating nodes}, which we call after the GMSH function \texttt{relocate\_nodes} that was used for this task. The strategy arises from the observation that the meshing step can be made continuous if the newly proposed $\vartheta^{\ast}$ is close to $\theta^{(l)}$. Concretely, we can use the \texttt{relocate\_nodes} function on the mesh $\mu(\vartheta^{(l)})$ to move its nodes around such that they fit the non-conducting region, described by $\vartheta^{\ast}$. The quality of the resulting mesh is checked using several quality measures built into GMSH. If these checks pass, this mesh is taken as the new mesh $\mu(\vartheta^{\ast})$. If they fail, a new independent mesh is generated for $\mu(\vartheta^{\ast})$ instead. In alignment with the superscript (IM) used for the discretized forward map based on independent meshing, we call this forward map $\hat{\mathcal{F}}^{(\text{RN},\vartheta=\vartheta^{(l)})}_{\Delta}(\vartheta)$.
 
\begin{figure}
    \centering
    \includegraphics[width=0.5\linewidth]{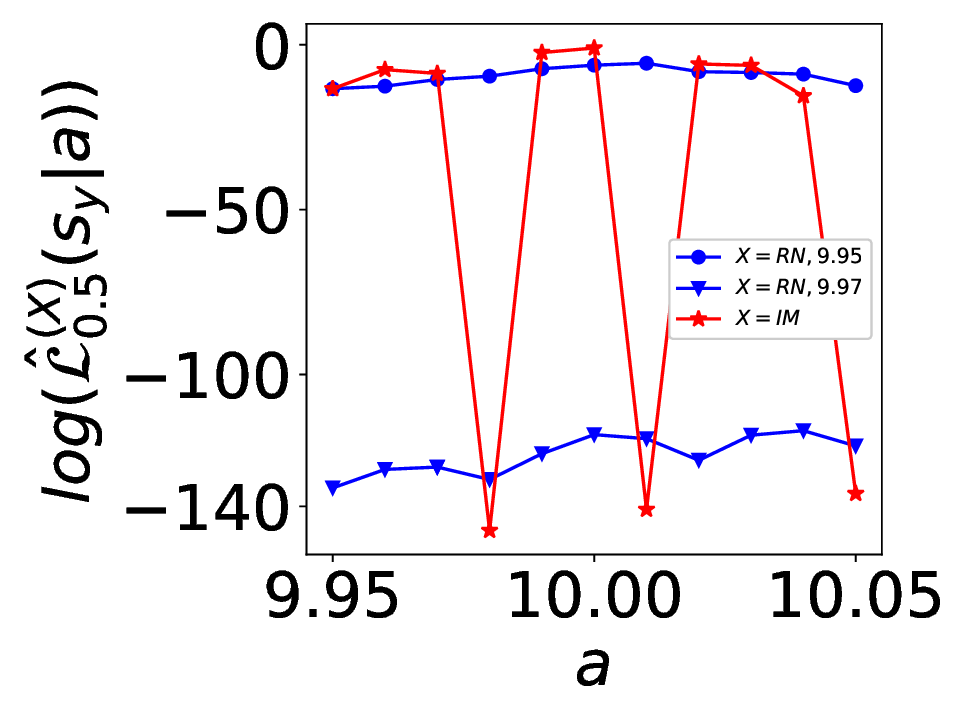}
    \caption{Log-likelihood curves generated with independent meshing $\hat{\mathcal{L}}_{0.5}^{(\text{IM})}(s_y|\vartheta)$ and Node Relocation, generated with $a=9.95$ and $a=9.97$ as starting point, $\hat{\mathcal{L}}_{0.5}^{(\text{RN},9.95)}(s_y|a)$ and $\hat{\mathcal{L}}_{0.5}^{(RN,9.97)}(s_y|\vartheta)$, respectively. The linear interpolations between the discrete evaluations of the likelihood are solely for visualization. }
    \label{fig:likelihood-ifo-meshing}
\end{figure}

 
We illustrate the benefit of this approach by showing what happens inside our adapted MH algorithm. We consider inference of the long radius, thus $\vartheta=a$, from the data $y$ from Section \ref{data}. We use the discretized compressed likelihood with $\Sigma_d=0$, hence, without modeling the term $d$ due to the discretization error. If the last obtained sample is $\vartheta^{(l)}=9.95$, Node Relocation offers a local continuous discretized compressed likelihood $\hat{\mathcal{L}}_{\Delta}^{(\text{RN},\vartheta=9.95)}(s_y|\vartheta)$, as can be seen in Figure \ref{fig:likelihood-ifo-meshing}. Meanwhile, if we use independent meshing, the discretized compressed likelihood $\hat{\mathcal{L}}_{\Delta}^{(\text{IM})}(s_y|\vartheta)$ strongly alternates solely due to the discretization error. The acceptance rate in (\ref{eq:ar-2}) depends on the relative difference in likelihood between $\vartheta^{(l)}$ and $\vartheta^{\ast}$, which means that its accuracy is strongly affected by such variations. For example, with independent meshing, proposing $\vartheta^{\ast}=9.97$ will almost certainly result in a rejection due to the large difference in the likelihood. Such rejections both affect the efficiency of the algorithm and the accuracy of the posterior. 

Finally, we consider the convergence properties of our adapted MH algorithm. AM is proven to be ergodic for the posterior distribution under minor assumptions, see Theorem 2 in \cite{haario2001}. Assumption (ii) of this theorem requires a fixed distribution to which can be converged. However, this is no longer the case. As we see in Figure \ref{fig:likelihood-ifo-meshing}, we obtain different discretized compressed likelihoods over the course of the algorithm. For example, if at another point in time, the adapted AM algorithm generates the sample $\vartheta^{(l)}=9.97$, then we obtain a discretized compressed likelihood $\hat{\mathcal{L}}_{\Delta}^{(\text{RN},\vartheta=9.97)}(\vartheta)$, locally continuous around $9.97$. The parameter value $\vartheta=9.95$ receives a very different likelihood value than for the likelihood around $\vartheta=9.95$. 

As a result, the posterior is no longer sampled exactly. However, as already mentioned, we do not have access to the exact posterior $P(\vartheta|y)$ anyway, but only to its discretized approximation $\hat{P}^{(IM)}_{\Delta}(\vartheta|y)$, which can be highly inaccurate if the discretization error is not modeled in the likelihood. Instead, our algorithm samples the posterior approximately. The benefit of our algorithm lies in the fact that it becomes more robust. If the likelihood accounts for the discretization error, the approximate samples closely match the exact samples but slightly overestimate the posterior variance. If the discretization error is not accounted for or underestimated, the approximate samples more closely resemble the non-discretized $P(\vartheta|y)$ than exact samples.

We show this using the discretized likelihood ratios. We define the pure data misfits, $\Delta y_{l}=y-\mathcal{F}(\vartheta^{(l)})$ and $\Delta y_{\ast}=y-\mathcal{F}(\vartheta^{\ast})$, which assume an exact model solution, and the discretization errors, $d_l=\mathcal{F}(\vartheta^{(l)})-\hat{\mathcal{F}}^{(X)}(\vartheta^{(l)})$ and $d_\ast=\mathcal{F}(\vartheta^{\ast})-\hat{\mathcal{F}}^{(X)}(\vartheta^{\ast})$,  for $\vartheta^{(l)}$ and $\vartheta{\ast}$, respectively. Similar as in \cite{lovbak2026}, we can rewrite the discretized likeliood ratio as as
\begin{align} 
   \log{ \frac{\hat{\mathcal{L}}^{(X)}_\Delta(s_y|\vartheta^{\ast})}{\hat{\mathcal{L}}^{(X)}_\Delta(s_y|\vartheta^{(l)})}}=& \Delta y_l^T\Sigma^{-1} \Delta y_l-\Delta y_{\ast}^T\Sigma^{-1} \Delta y_{\ast} \nonumber\\&-2 \Delta y_l^T\Sigma^{-1} d_l+2\Delta y_\ast^T\Sigma^{-1} d_{\ast}-d_l\Sigma^{-1} d_l +d_{\ast}\Sigma^{-1} d_{\ast} \label{eq:ar-discr}\\
    =& \log{ \frac{\mathcal{L}(s_y|\vartheta^{\ast})}{\mathcal{L}(s_y|\vartheta^{(l)})}} -2 \Delta y_l^T\Sigma^{-1} d_l+2\Delta y_\ast^T\Sigma^{-1} d_{\ast}-d_l\Sigma^{-1} d_l +d_{\ast}\Sigma^{-1} d_{\ast}
\end{align}
The superscript $(X)$ refers to the discretization technique. First, we note that for Node Relocation $\lim_{\vartheta^{\ast}-\vartheta^{(l)}\rightarrow 0} \log(\hat{a}^{(X)}_\Delta(\vartheta^{(l)},\vartheta^{\ast})$ exists and equals 0. In general, we can distinguish two possible scenarios during the accept-reject step:
\begin{itemize}
    \item \textbf{Case 1: $\mathbf{|\theta^{\ast} - \theta^{(l)}|}$ is small.} We expect a similar data misfit for both parameter values. As a result, the relative impact of $d_l$ and $d_\ast$ on (\ref{eq:ar-discr} is large. Node Relocation now reduces this impact and guarantees that $d_{\ast}\rightarrow d_l$ for $\vartheta^{\ast}\rightarrow \vartheta^{(l)}$. This is especially desired if the likelihood does not model the discretization error, because then it closely matches the exact acceptance rate. If the likelihood does not model the discretization error, correlating $d_l$ and $d_\ast$ actually causes the variance of the discretization error to be slightly too small in comparison to the likelihood model, such that we might slightly overestimate posterior variance.
    \item \textbf{Case 2: $\mathbf{|\theta^{\ast} - \theta^{(l)}|}$ is large.} We cannot generate the mesh with Node Relocation such that both discretization techniques result in the same acceptance probability. However, $\Delta y_l$ and $\Delta y_\ast$ differ strongly, such that the first two terms dominate in (\ref{eq:ar-discr})
\end{itemize}

\section{Computational results}
\todo[inline]{Rerun experiments with improved conditions}
We demonstrate our Bayesian inference setup with the discretized compressed likelihood and our adapted MH algorithm on the synthetic intra-atrial reentrant tachycardia (IART) test case. First, we dedicate Section \ref{results-exp-1} to evaluating our algorithm from Section \ref{MCMC}. We test the Adaptive Metropolis proposal and the Node Relocation meshing strategy in the accept-reject step of our adapted MH algorithm by comparing it to RWMH. In Section \ref{results-exp-2}, we solve the inference experiment defined in Section \ref{Bayes} and estimate the geometrical parameter from the data to illustrate Bayesian inference in this context. In Section \ref{results-exp-3}, we repeat this experiment for a compressed likelihood based on a more conservative noise model. Due to the increased amount of noise, the inference problem suffers from identifiability issues.  However, Bayesian inference still performs accurately in this scenario and offers additional insights into the identifiability issues.

\subsection{Experiment 1: Evaluation and performance of the adapted MH algorithm}\label{results-exp-1}
Section \ref{MCMC-proposal} details the proposal step of our adapted MH algorithm. We use an Adaptive Metropolis proposal with a covariance equal to the sample covariance of the obtained samples. Given that the joint posterior distribution of $a$ and $b$ is highly correlated, adaptation of the proposal covariance to match the posterior covariance structure can strongly improve the sampling efficiency. We illustrate this for the inference experiment laid out in Section \ref{Bayes} with $\gamma=0.8$ and $dx=0.25$~mm and compare our adapted MH to RWMH. We take the same manually tuned $\Sigma_0=\diag(0.05,0.05,0.01)$ for both. For our method, the covariance will change iteratively according to Equation (\ref{eq:am}), while it will remain constant for RWMH.
    
\begin{figure}[ht]
    \centering
      \includegraphics[width=0.49\linewidth]{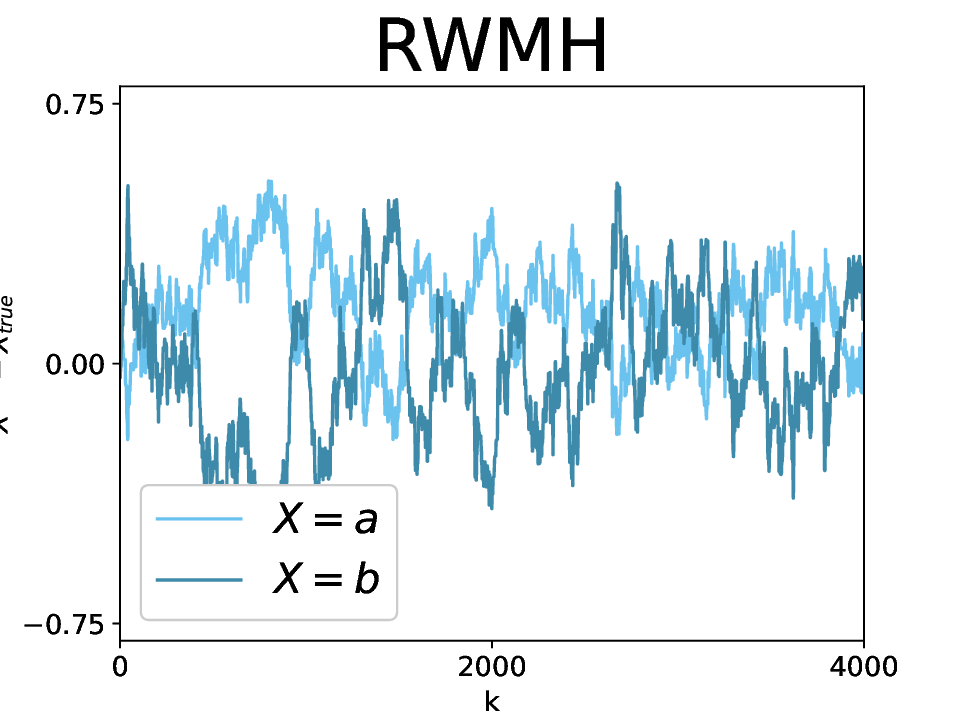}        
      \includegraphics[width=0.49\linewidth]{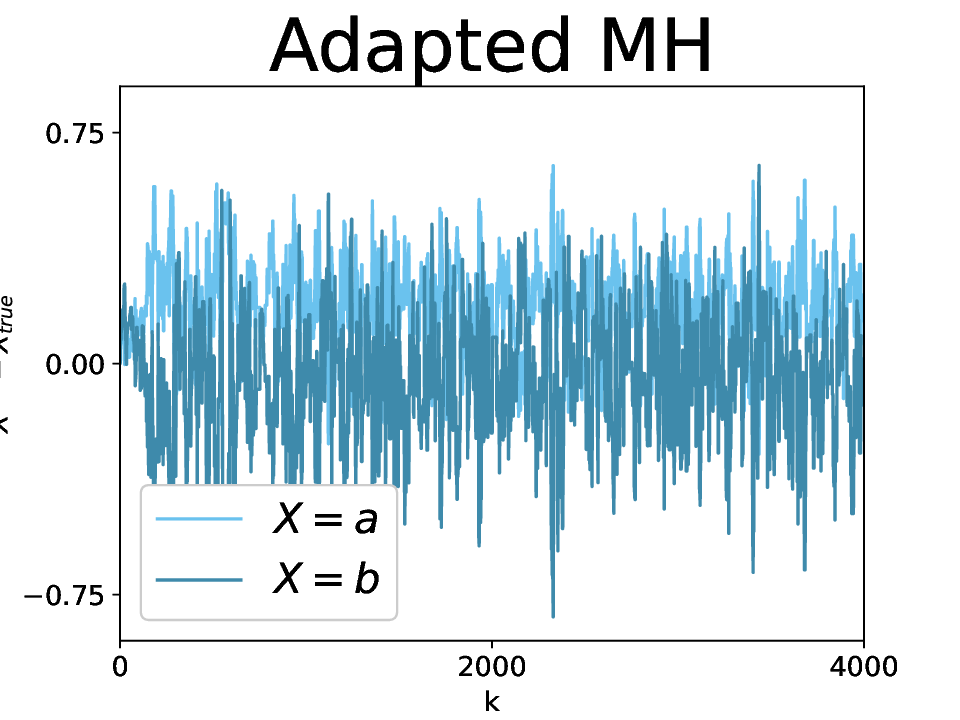}    
    \caption{The $a$- and $b$-component of the Markov chain $(\vartheta-\vartheta_{true})$ inside the MCMC method. Left: bad mixing with RWMH; right: good mixing with AM.}
    \label{fig:mixing}
\end{figure}

Figure \ref{fig:mixing} illustrates the superior mixing of the Markov chains, and, thus, the superior parameter space exploration using our method in comparison with RWMH. Thanks to the adaptive covariance structure, the principal axes of the proposal align with the principal axes of the posterior such that larger moves are proposed in directions with larger variance. In contrast, RWMH overestimates the posterior variance in certain directions while underestimating it in others. To achieve a good acceptance rate, small steps are taken in directions with more variance to compensate for taking large steps in directions with less variance. The former are almost always accepted, while the latter are rejected frequently. Figure \ref{fig:mixing} illustrates how the resulting chains mix poorly. The chains cannot sample the entire posterior continuously, but instead sample regions of the posterior sequentially. The consequence of this behavior is that subsequent samples are correlated, and, as said in Section \ref{MCMC-proposal}, this high autocorrelation results in slower convergence of the resulting estimations. 
 
Next, we compare Node Relocation and Independent Meshing. For the forward model, we choose $\gamma=0.1$. We now only infer $a$, thus $\vartheta=a$ in all expressions in Sections \ref{Bayes} and \ref{MCMC}. The proposal distribution $f_u$ in Equation (\ref{eq:proposal}) simplifies to $f_u\sim\mathcal{N}(0,\sigma_f^2)$ and we tune $\sigma_f$ separately for Experiments 1a, 1b and 1c (see Table \ref{tab:acceptance-rates}) using the same value for Independent Meshing and Node Relocation. For each experiment we set the number of samples to $N=7\,500$. In Experiment 1a, we take $\Delta x=0.25$~mm and use a correct estimate for the discretization error in the discretized compressed likelihood with the procedure defined in Section \ref{Covariance}. In Experiments 1b and 1c, we coarsen the resolution to $\Delta x=0.5$ without adapting $\Sigma_d$, such that the discretization error is underestimated. 

\begin{figure}[ht]
    \centering
    \includegraphics[width=\linewidth]{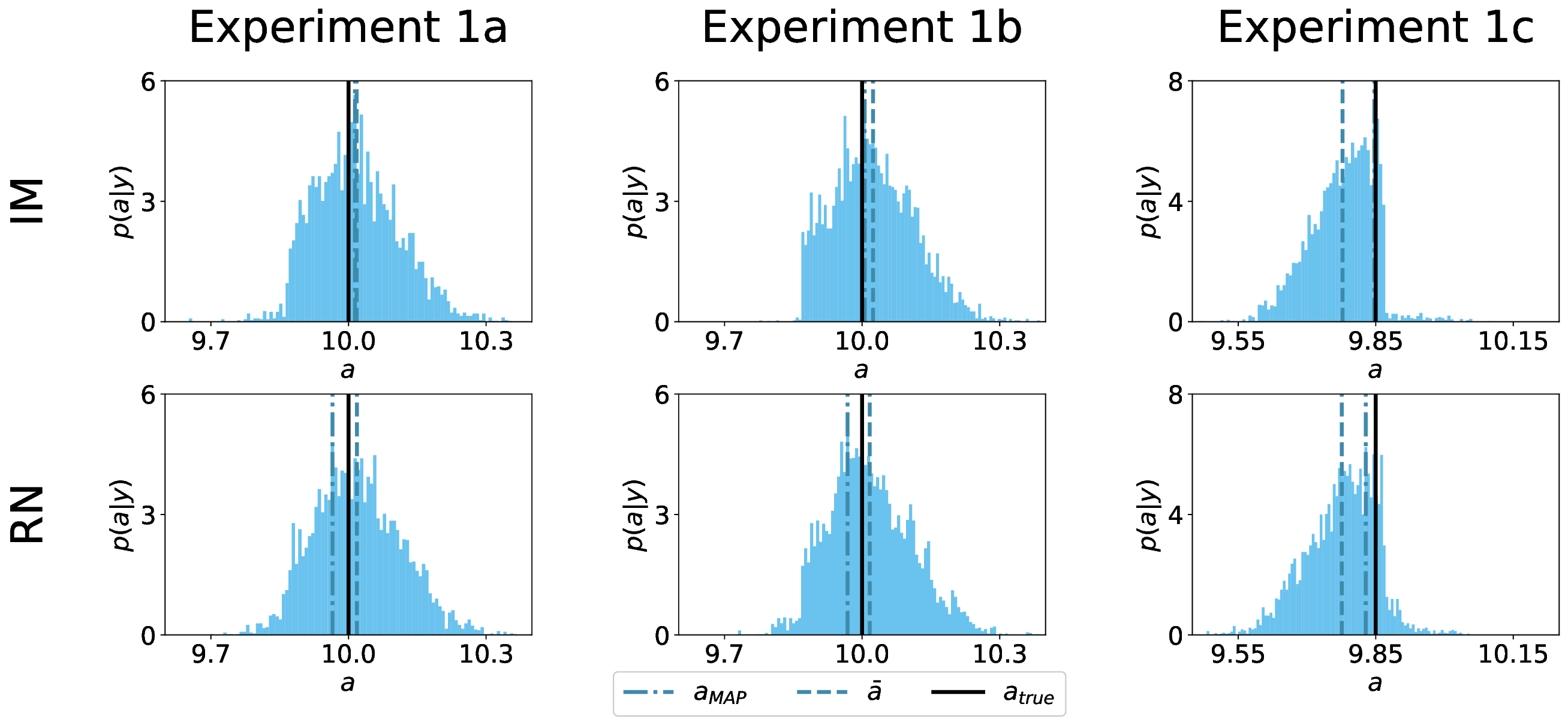} 

    \caption{Inference of $a$ with $a_{true}=10$ mm using Node Relocation (RN) and Independent Meshing (IM) in the accept-reject step. Experiment 1a: $dx=0.25$, $\Sigma_d$ modeled as proposed in Section \ref{Covariance}; Experiment 1b: same covariance but $dx=0.5$; Experiment 1c: same covariance, but $dx=0.5$ and $a_{true}=9.85$~mm}    
    \label{fig:1c}
\end{figure}

\begin{figure}[ht]
     \centering
    \includegraphics[width=0.45\linewidth]{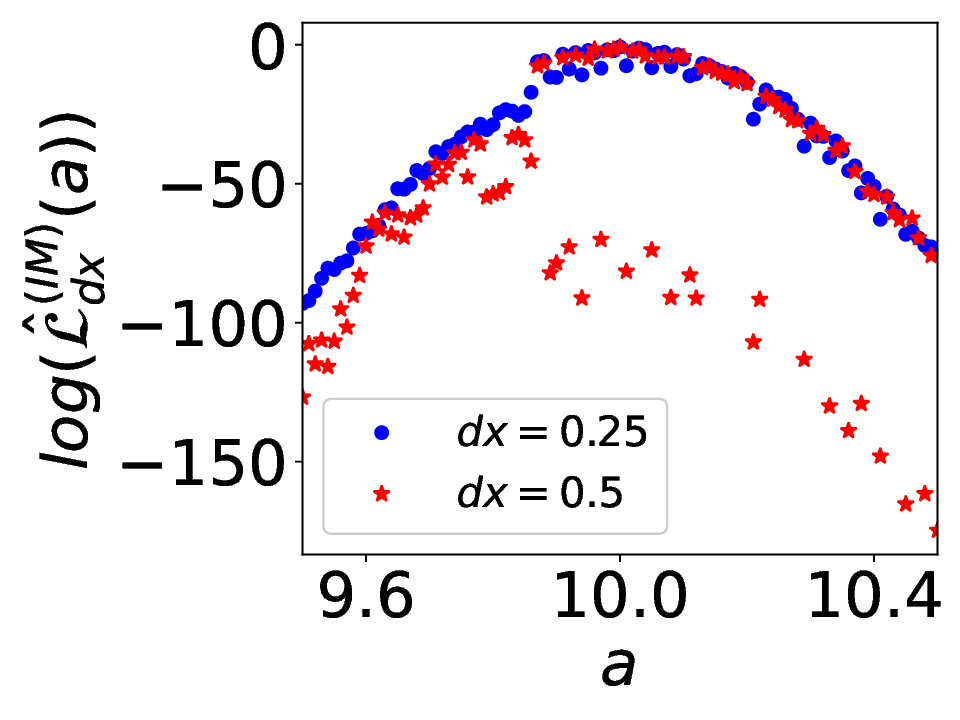}
     \caption{The discretized compressed log likelihood for the same inference problem as in Figure \ref{disc-ll}. Now the variance is $\Sigma_{\varepsilon'}+\Sigma_d$, with $\Sigma_d$ estimated for $dx=0.25$.}
     \label{fig:ll-plot-with-sd}
 \end{figure}

Figure \ref{fig:ll-plot-with-sd} shows the discretized compressed likelihood for this covariance. In this plot, we see two types of discontinuities. On the one hand, the discretized compressed likelihood visually 'oscillates' due to the discretization error over part of the parameter range, alternating between a higher located curve and a lower located one. On the other hand, an abrupt discontinuity appears around $a=9.85$~mm. We choose Experiment 1b and 1c with $a_{\text{true}}=10.0$ and $a_{\text{true}}=9.85$, such that this abrupt discontinuity occurs once in the tail and once close to the mode of the posterior distribution, respectively. The resulting posteriors and mean estimates are shown in Figure \ref{fig:1c} and in Table \ref{tab:acceptance-rates}.

As we explained in Section \ref{MCMC-ar}, the acceptance probability is highly inaccurate when using Independent Meshing if we underestimate the discretization error. This can strongly affect the accuracy of the posterior histogram. The most visible effect is noticed at the point where the discretized compressed likelihood has an abrupt cut-off for the coarse discretization. Experiment 1b and 1c show how this abrupt cut-off is translated into the posterior. For Independent Meshing, we observe an abrupt cut-off to almost zero probability around $a=9.85$ as well. Meanwhile, for Node Relocation, this cut-off is somewhat smoothed such that the posterior shape is modeled slightly better. However, even if the discontinuity is smoothed, it does not change the fact that the discretization error is underestimated. At least, Independent Meshing reveals visually when this error is underestimated and when better modeling is required.

The `oscillating" discontinuities will result in similar inaccuracies when comparing one sample from the top of an oscillation with one from the bottom. While these are not directly visible, they will still affect the histogram shape. After all, the histogram calculates locally which fraction of the posterior samples occurs within a bin. As the vast majority of accepted samples have a discretized compressed likelihood value on the upper end of the oscillations, these fractions will be affected by the fraction of the parameter values within the bin width that has a large discretized compressed likelihood value. As this can change over the parameter space, this can affect the accuracy of the histogram. Despite this effect, we observe that the mean estimates of Node Relocation and Independent Meshing are indistinguishable.

The accuracy of the acceptance probability also affects the efficiency. Independent Meshing can cause a large number of rejections based on large differences in the forward solution that are mainly caused by the discretization error. Meanwhile, Node Relocation avoids many of these futile forward simulations. Table \ref{tab:acceptance-rates} confirms this, as the acceptance rates are higher for the same proposal distribution for Node Relocation.

\begin{table}[ht]
    \centering
     \caption{Acceptance rates and mean estimates for both meshing strategies for the three experiments. IM = Independent Meshing; NR = Node Relocation}
    \begin{tabular}{lccccc}
        \toprule
        & $\mathbf{\sigma_f}$& \multicolumn{2}{c}{\textbf{Acceptance Rate}} & \multicolumn{2}{c}{\textbf{Mean estimate}} \\
        \cmidrule(lr){3-4} \cmidrule(lr){5-6}
        & &\textbf{IM} & \textbf{NR}  & \textbf{IM} & \textbf{NR} \\
        \midrule
           \textbf{Experiment 1a}& 0.2&0.440 & 0.452 & 10.02 & 10.02 \\
        \textbf{Experiment 1b}&0.15&0.411 & 0.440& 10.02 & 10.02 \\

        \textbf{Experiment 1c}&0.15&0.402 & 0.413 & 9.78 & 9.78 \\

        \bottomrule
    \end{tabular}
   
    \label{tab:acceptance-rates}
\end{table}

We conclude that correct discretization error estimation is more important than the chosen meshing strategy. If it is possible to visually check the posterior shape, Independent Meshing might be the better choice because it will warn when the discretization error is underestimated. However, for higher-dimensional parameters, visual checks become more complicated, and we could opt for Node Relocation because of its higher efficiency and robustness.
\subsection{Experiment 2: inference based on one catheter location} \label{results-exp-2}
We solve the Bayesian formulation of our parameter inference problem, defined in Section \ref{Bayes}. We use a discretization $dx=0.25$ and take the conduction slowing $\gamma=0.8$. For these values, the discretization error can be ignored relative to the measurement noise ($\Sigma = \Sigma_{\varepsilon'}$). We determine the posterior using our adapted MH algorithm, with an adaptive proposal and a Node Relocation strategy in the accept-reject step. We choose $\Sigma_0=\diag(0.0025,0.0025,0.0001)$ in (\ref{eq:am}) and $N=10\,000$.

\begin{figure}[ht]
    \includegraphics[width=0.32\linewidth]{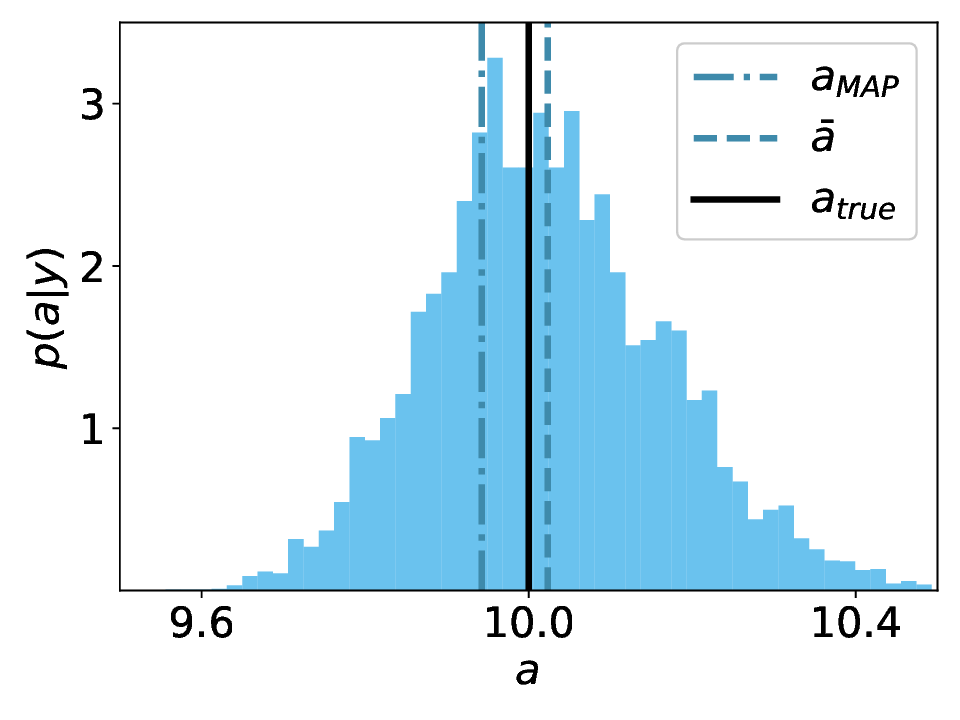}        
      \includegraphics[width=0.32\linewidth]{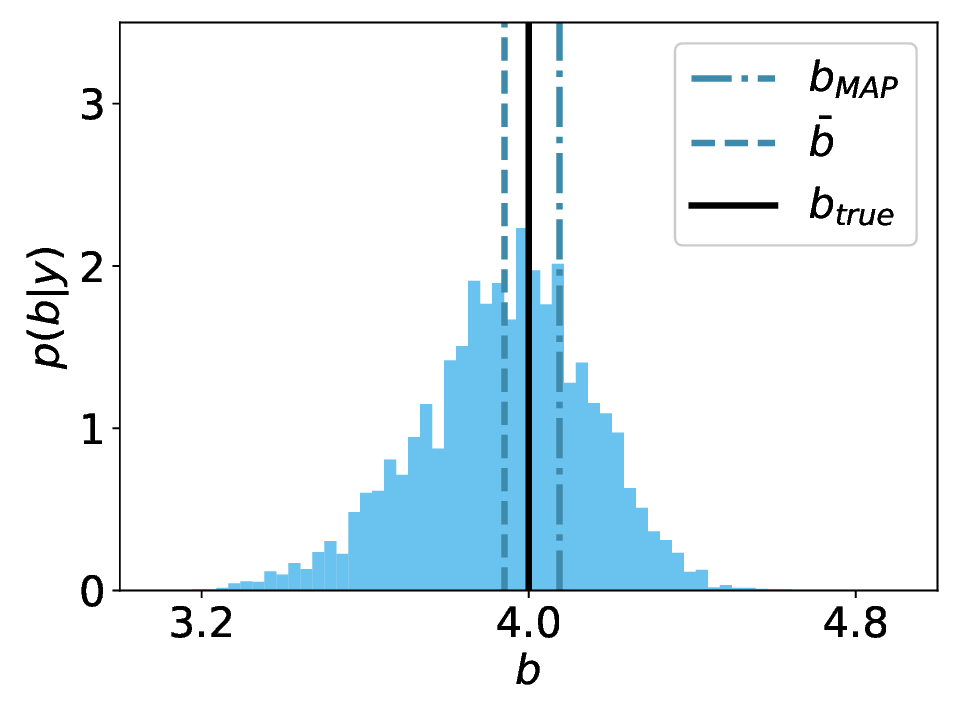}             
    \includegraphics[width=0.32\linewidth]{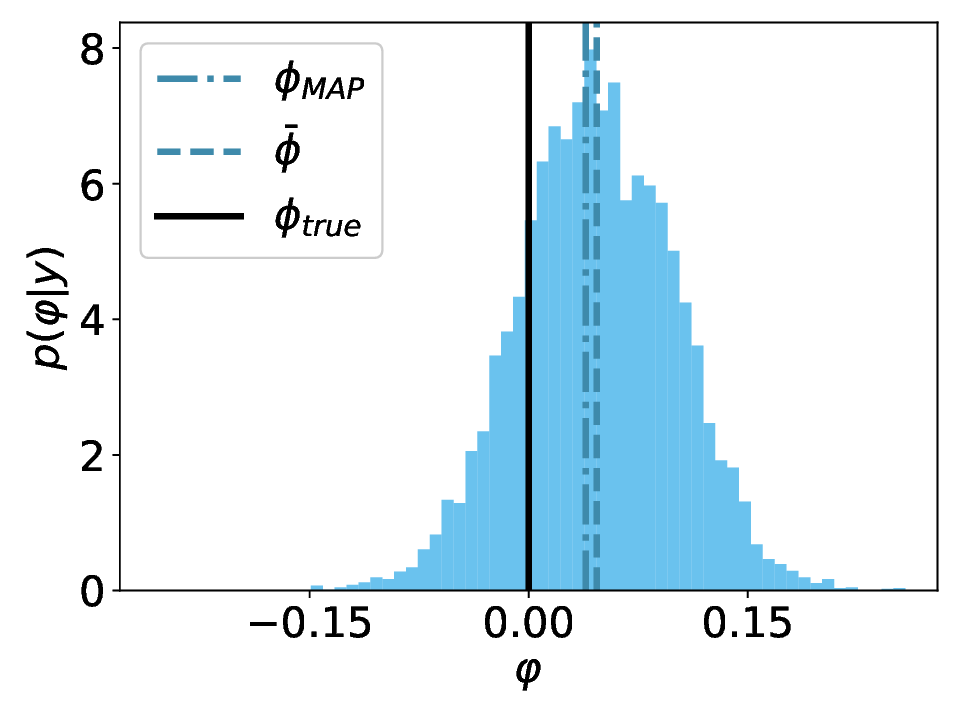}             
    \caption{Histograms that approximate the marginal posteriors for inference of $\vartheta$ based on the EGM data measured at one catheter location (Experiment 2). Left: $p(a|y)$; middle: $p(b|y)$; right: $p(\varphi|y)$. The dashed line gives the mean estimate, and the dashdotted line the Maximum A Posteriori (MAP) estimate for the parameter value.}
    \label{fig:1-detector-a}
\end{figure}

Figure \ref{fig:1-detector-a} shows the histograms that approximate the marginal posterior on $a$, $b$ and $\varphi$. Table \ref{tab:estimations-exp3-1} contains the mean estimate, the Maximum A Posteriori (MAP) estimate and the sample variances that reflect the uncertainty on the individual parameters $a$, $b$ and $\varphi$. By using Bayesian inference, we obtain both an estimate, namely either the mean or the MAP, and an uncertainty on this estimate, namely the variance.  While the measurement error causes an error in the estimates, the true value still has a large posterior probability.
 
\subsection{Experiment 3: Bayesian inference in presence of identifiability issues}\label{results-exp-3}
In the last section, we repeat the same experiment as in Section \ref{results-exp-2}, but now we assume an increased amount of measurement noise on the characterizing quantities in the compressed discretized likelihood. For this purpose, we set $\Sigma=\Sigma_{\varepsilon'}=10.0I$ and $\Sigma=\Sigma_{\varepsilon'}=\diag(2.5,25.0,\dots,25.0)$ in Experiment 3a and 3b, respectively. For Experiment 3c, we also reduce the diffusion by setting $\gamma=0.1$, while using a covariance matrix $\Sigma=10.0I$, which now accounts for measurement error and discretization error. The relevance of these experiments stems from the fact that we are deriving a more realistic noise model from clinical data to be used in follow-up work. This work might reveal that we underestimated the measurement noise considerably. To perform sampling, we set $\Sigma_0=\diag(0.25,0.25,0.05)$ for Experiment 3a and 3c and $\Sigma_0=\diag(0.8,0.8,0.16)$ for Experiment 3b. We set $N=10\,000$ for 3a and 3b, and $N=7\,500$ for 3c.
 
\begin{figure}[ht]
    \centering
     \includegraphics[width=0.32\linewidth]{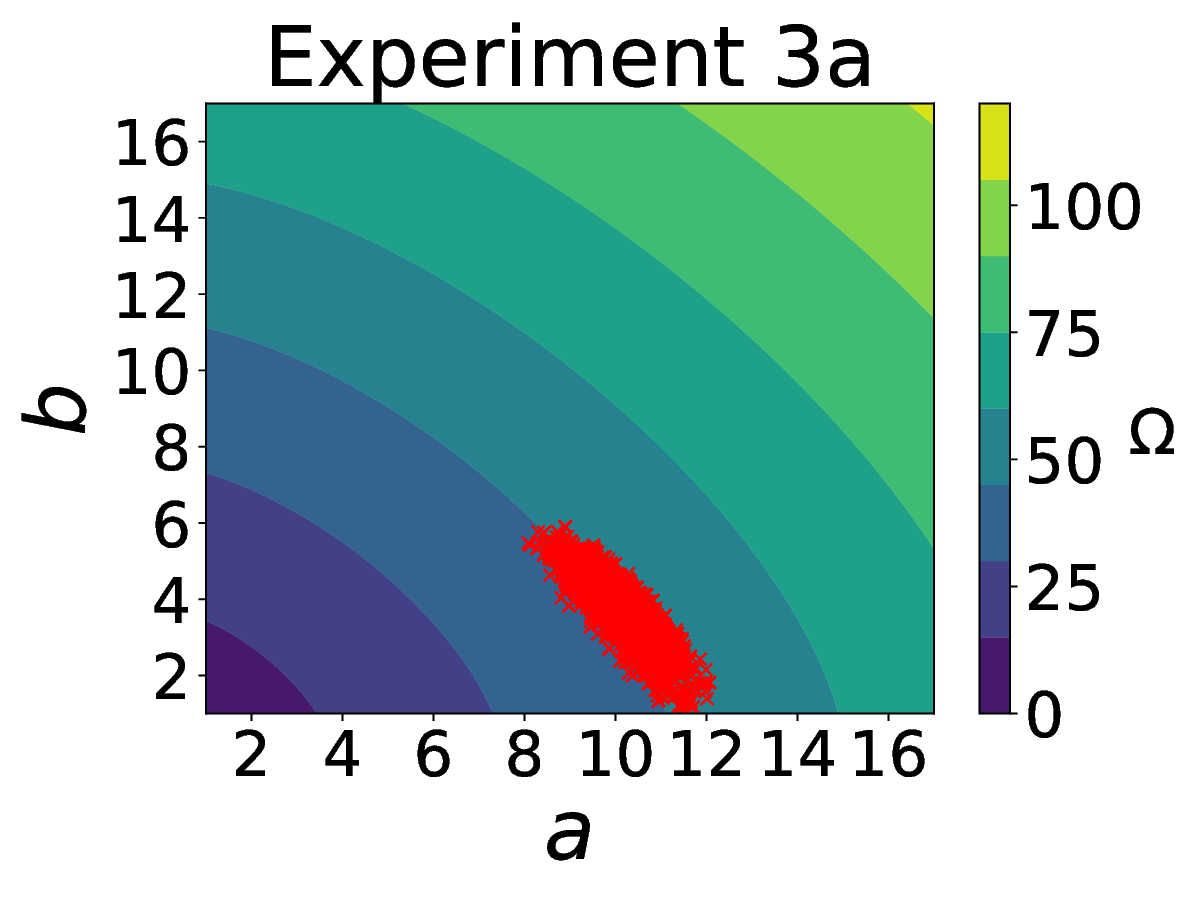}        
     \includegraphics[width=0.32\linewidth]{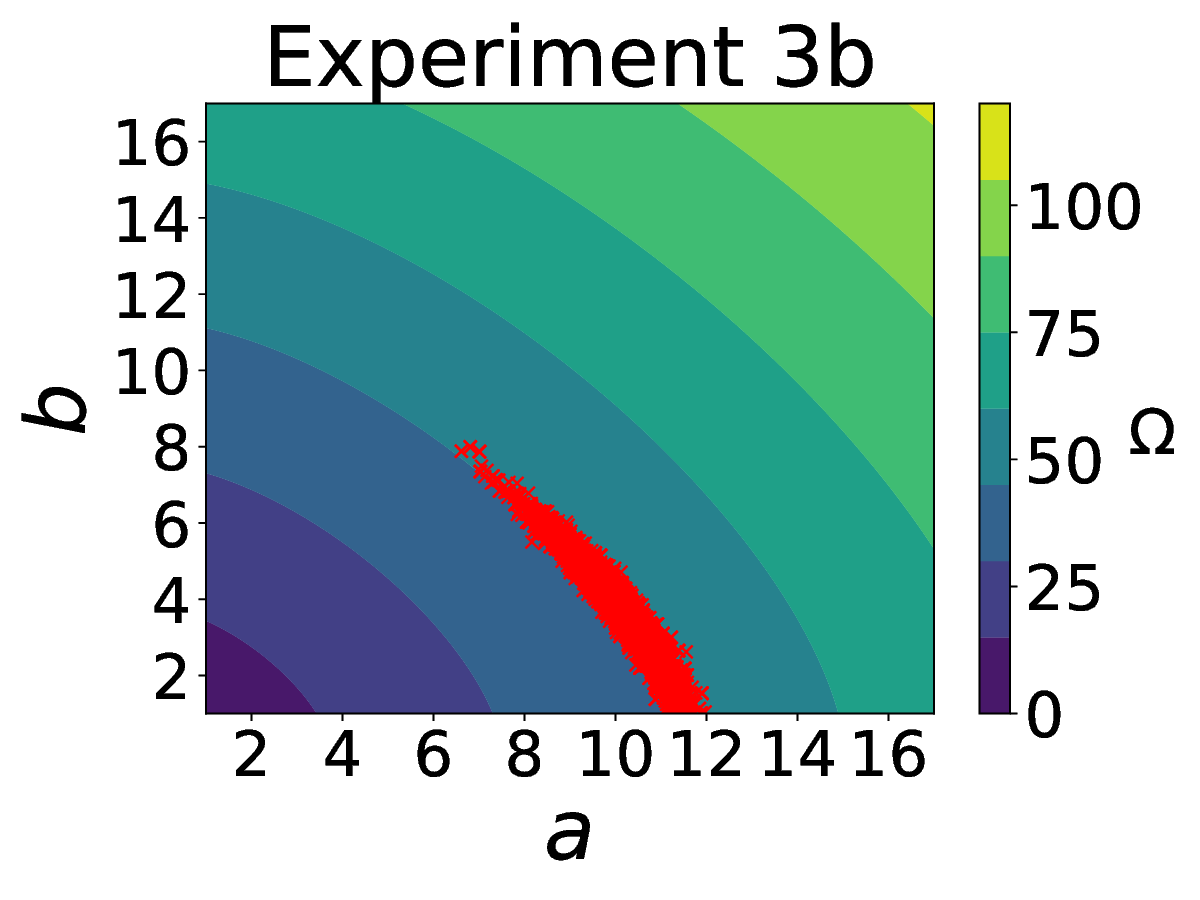}        
     \includegraphics[width=0.32\linewidth]{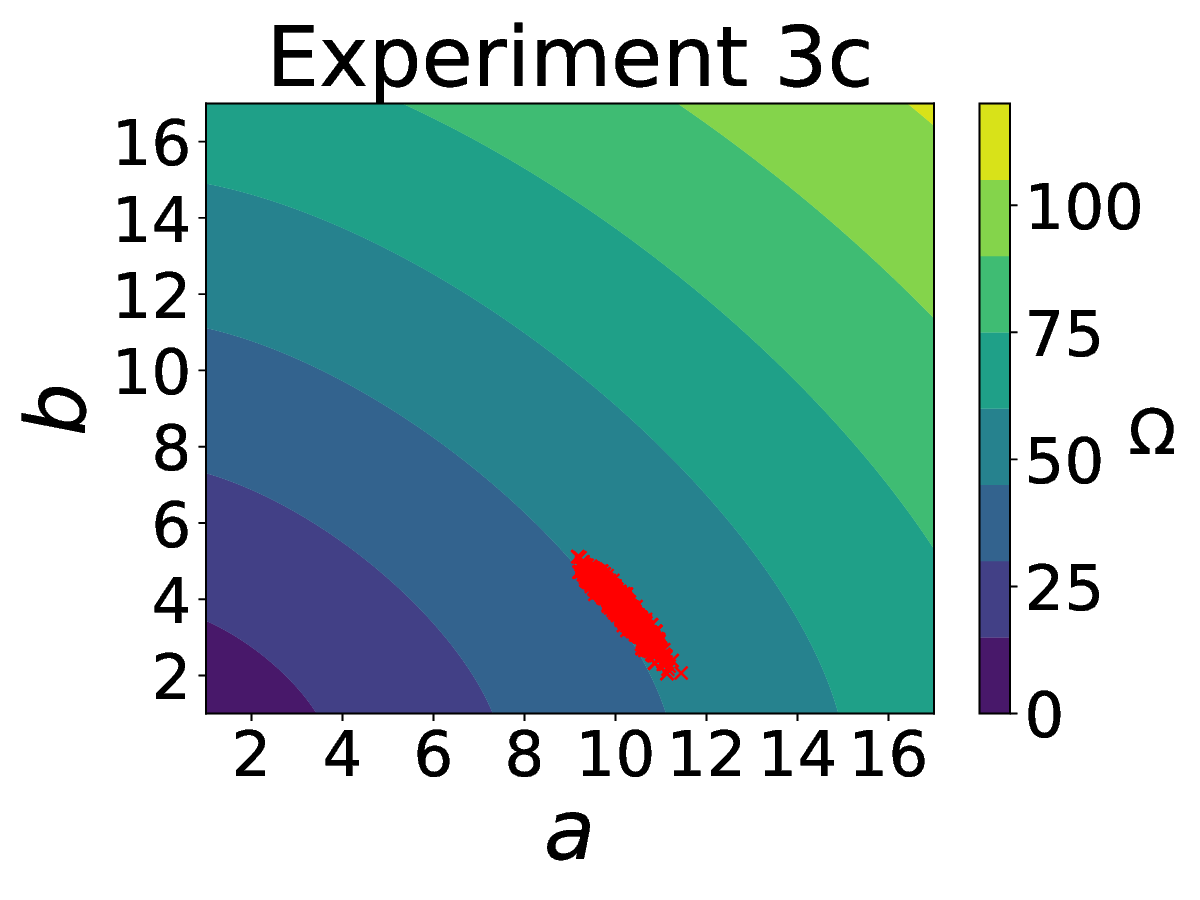}        

    \caption{Samples of the posterior distribution (black crosses) for Experiment 3 in the ab-plane, together with perimeter isolines.}
    \label{fig:exp3-1a}
\end{figure}

 Using Bayesian inference and our adapted MH algorithm, we obtain the posterior distribution on $\vartheta$, similarly to before. We do not show all the histograms, but plot the obtained samples in the $(a,b)$-plane in Figure \ref{fig:exp3-1a}. Moreover, Table \ref{tab:estimations-exp3-1} shows the resulting parameter estimates and the variances on $a$, $b$ and $c$. From these results, we observe that a large amount of uncertainty remains on $a$ and $b$, and we conclude that we cannot infer much about $a$ and $b$ from the information in the data. We can make some estimates, but these have a bias, with much uncertainty remaining on this estimate. However, this is considerably more accurate than making a wrong estimate on its own as in deterministic approaches. 
 
 Furthermore, we observe that nearly all uncertainty in $a$ and $b$ lies along a single direction in the $(a,b)$-plane (see Figure \ref{fig:exp3-1a}). Starting from a uniform, uncorrelated prior, the data reduces uncertainty to almost one dimension, aligned with the direction of constant perimeter. Consequently, all high-probability combinations of $a$ and $b$ yield similar perimeters. Because this direction does not align with either axis, substantial uncertainty remains for each parameter individually, though the joint uncertainty is less than the sum of their marginal uncertainties.

\begin{table}[ht]
    \centering
        \caption{Mean and MAP estimate and  the variance in the direction of the ellipse axes $a$ and $b$ and the inclination angle $\varphi$ for Experiment 2 and 3.}
    \begin{tabular}{lccccc}
        \toprule
     &$\Bar{\vartheta}$&$\vartheta_{MAP}$&$S_a^2$ & $S_b^2$ &$S_{\varphi}^2$  \\
        \midrule
        \textbf{Exp. 2}&$\begin{bmatrix}
        10.02\; \text{mm}\\
        3.94\; \text{mm}\\
        0.05\; \text{rad}
    \end{bmatrix}$ &$\begin{bmatrix}
        9.94 \; \text{mm}\\ 4.08\; \text{mm}\\0.04 \;\text{rad}
    \end{bmatrix}$&2.05e-02 & 4.45e-02 & 2.85e-03\\
    \midrule
        \textbf{Exp. 3a}
        &$\begin{bmatrix}
        10.12\; \text{mm}\\
        3.76\; \text{mm}\\
        0.07\; \text{rad}
    \end{bmatrix}$&$\begin{bmatrix}
           9.63\; \text{mm}\\
        4.26\; \text{mm}\\
        -0.01\; \text{rad}
    \end{bmatrix}$&3.23e-01 & 5.51e-01 & 2.50e-02\\
    \midrule
    \textbf{Exp. 3b}
        &$\begin{bmatrix}
        10.08\; \text{mm}\\
        3.76\; \text{mm}\\
        0.09\; \text{rad}
    \end{bmatrix}$&$\begin{bmatrix}
        10.16\; \text{mm}\\
        3.89\; \text{mm}\\
        0.00\; \text{rad}
    \end{bmatrix}$&5.58e-01 & 1.31e+00 & 5.17e-02\\
        \midrule
    \textbf{Exp. 3c}
        &$\begin{bmatrix}
        10.13\; \text{mm}\\
        3.81\; \text{mm}\\
        0.05\; \text{rad}
    \end{bmatrix}$&$\begin{bmatrix}
        10.11\; \text{mm}\\
        3.85\; \text{mm}\\
        0.06\; \text{rad}
    \end{bmatrix}$&1.17e-01 & 2.16e-01 & 7.31e-03 \\
    \bottomrule
    \end{tabular}

    \label{tab:estimations-exp3-1}
    \end{table}

 These results nicely illustrate the benefits of the stochastic fits over deterministic fits. We observe that a deterministic estimate is undeniably worse than the estimates in Section \ref{results-exp-2}. We could falsely assume that we fully characterized the non-conducting region, even if we would not really know $a$ and $b$. Another benefit of Bayesian inference is that it can offer more insights than a deterministic fit. For example, in our setup, the Bayesian inference results reveal identifiability issues for $a$ and $b$ from the data. At the same time, we see that we can already estimate the perimeter of the non-conducting region, $\Omega$, which could have gone unnoticed in a deterministic approach. 

\section{Conclusions and outlook} This paper demonstrated Bayesian inference of a stochastic non-conducting region parametrization from synthetic electrogram (EGM) data. We illustrated the benefits of this approach compared to deterministic fitting. On top of an estimate, it quantifies the uncertainty in the inferred geometric parameters, providing a measure of reliability. Additional insights are obtained, such as identifiability issues. 

We proposed a methodology to minimize the number of samples required for inference by using an adapted MH algorithm for MCMC sampling. The latter combines an Adaptive Metropolis proposal strategy and an accept-reject step using node relocation for 
robust and efficient sampling under posterior correlations and discretization error. However, the cost per sample remains high, as an expensive forward solve of the monodomain equation is required for each sample. As the number of required samples increases with the parameter dimension, inference of high-dimensional parameter sets may become computationally prohibitive. Future work will focus on strategies that replace many of these expensive simulations with computationally cheaper approximations. For this line of research, the presented work can act as the ideal starting point.

This study also examined the impact of the uncertainty induced by discretization error. We analyzed the behavior of the forward model and the posterior under discretization error during inference of the geometrical parameter. Furthermore, we compared Independent Meshing and Node Relocation within the accept–reject step of the MCMC algorithm. While node relocation offers a more robust and faster estimation, independent meshing can visually reveal underestimated discretization errors and the associated overconfidence in predictions.

We defined a likelihood based on characterizing quantities for improved accuracy of the inference. One limitation concerns the exhaustive modeling of measurement uncertainty. This primarily relates to the noise model of the characteristic-quantity vector, which serves as input to the Bayesian inference framework. One motivation for adopting the compressed likelihood was to facilitate such a noise model. The proposed methodology is directly applicable once this model is specified, though the posterior distribution may change substantially. An accurate noise model will enable the application of the framework to experimental data and a detailed analysis of the uncertainties introduced by measurement and modeling procedures.

\section*{Acknowledgments}
Our work was partially funded by the KU Leuven Research Council through grant IDN/22/006 and the Research Foundation – Flanders (FWO) through grants 1S53625N and 11PMS24N.

\appendix
\section{Coordinates of the electrodes of the mapping catheter at Location 1 and 2} \label{appendix-loc}
\vspace{1cm}
    \begin{center}
\begin{tabular}{lcc}
    \toprule
    Number & $x$-coordinate & $y$-coordinate \\
    \midrule
    1 & 27.50 & 62.00 \\
    2 & 32.50 & 62.00 \\
    3 & 37.50 & 62.00 \\
    4 & 42.50 & 62.00 \\
    5 & 25.77 & 64.38 \\
    6 & 27.32 & 69.13 \\
    7 & 28.86 & 73.89 \\
    8 & 30.41 & 78.64 \\
    9 & 22.98 & 63.47 \\
    10 & 18.93 & 66.41 \\
    11 & 14.89 & 69.35 \\
    12 & 10.84 & 72.29 \\
    13 & 22.98 & 60.53 \\
    14 & 18.93 & 57.59 \\
    15 & 14.89 & 54.65 \\
    16 & 10.84 & 51.71 \\
    17 & 25.77 & 59.62 \\
    18 & 27.32 & 54.87 \\
    19 & 28.86 & 50.11 \\
    20 & 30.41 & 45.36 \\
    \bottomrule
\end{tabular}

     \end{center}

\section{Calculation of Local Activation Times (LAT)}

To calculate the relLATs and the period, we first need to derive the local activation time (LAT) at each electrode, which is the arrival time of the spiral wave front at that electrode. For this step, we distinguish two situations. We illustrate this in Figure \ref{fig:LAT}
\begin{itemize}
    \item \textbf{Calculation of LATs from electrogram data $\boldsymbol{y}$}. Let us call the LAT at electrode $j$ for the first and second reentry $LAT1^{y}_{j}$ and $LAT2^{y}_{j}$, respectively. The electrogram time trace at that electrode is given in row $j$ of the matrix $y$ in Equation (\ref{eq:y}):
      \begin{equation}
    [y(r_{j},\tau_0),y(r_{j},\tau_1),\dots,y(r_{j},\tau_{N_{\text{frame}}})].
    \end{equation}
    The LATs are the zero-crossings with negative descent of this time trace. We find these by looking for the index $i=m$ in this time trace such that $y(r_{j},\tau_{m-2:m})>0$ and $y(r_{j},\tau_{m+1:m+3})<0$. To ignore zero-crossings caused by large noise realizations, we require that multiple measurements before $\tau_m$ are positive and multiple measurements after $\tau_m$ are negative. For the smallest value of $m$, we calculate $LAT1^{y}_j$ using linear interpolation as 
    \begin{equation}
        LAT1^{y}_j=\tau_m+\frac{y(r_{j},\tau_m)}{y(r_{j},\tau_m)-y(r_{j},\tau_{m+1})}\Delta \tau
    \end{equation}
    with $\Delta \tau$ the time between measurements. An equivalent formulation holds for $LAT2^{y}_j$ using the second value for $m$. 
    \item \textbf{Calculation of LATs from the discretized solution $\boldsymbol{\mathcal{F}(\vartheta)}$}. Let us call the prediction by the model for the LAT at electrode $j$ during the first and second reentry $LAT1^{\mathcal{F}(\vartheta)}_j$ and $LAT2^{\mathcal{F}(\vartheta)}_j$, respectively. Since $V_m(x,t)$ is only defined in the tissue and not at the electrode location outside the tissue, we have to project the location of electrode $j$ orthogonally onto the 2D-plane, $x_j=r_{j,1:2}$. From $\hat{\mathcal{F}}(\vartheta)=[\hat{V}_m(x,t_0,\vartheta),\hat{V}_m(x,t_1,\vartheta),\dots,\hat{V}_m(x,t_N,\vartheta)]$, we can now compose a transmembrane potential time trace at $x_j$ as
    \begin{equation}
    [\hat{V}_m(x_{j},\tau_0,\vartheta),\hat{V}_m(x_{j},\tau_1,\vartheta),\dots,\hat{V}_m(x_{j},\tau_{N_{\text{frame}}},\vartheta)].
    \end{equation}
    To obtain the LATs, we find the index $i=m$ such that the transmembrane potential exceeds 0.3 between $\tau_m$ and $\tau_{m+1}$, thus $\hat{V}_m(x_{j},\tau_m,\vartheta)<0.3$ and $\hat{V}_m(x_{j},\tau_{m+1},\vartheta)>0.3$. We use a similar interpolation, which, for the lowest value of $m$, gives us the expression
    \begin{equation}
        LAT1_j^{\mathcal{F}(\vartheta)}=\tau_m+\frac{0.3-V_m(x_{j},\tau_m,\vartheta)}{V_m(x_{j},\tau_{m+1},\vartheta)-V_m(r_{k,j},\tau_{m},\vartheta)} \Delta\tau,
        \end{equation}
        with a similar interpolation used for $LAT2_j^{\hat{\mathcal{F}}_\Delta(\vartheta)}$.

 \begin{figure}[ht]
     \centering
    \includegraphics[width=0.49\linewidth]{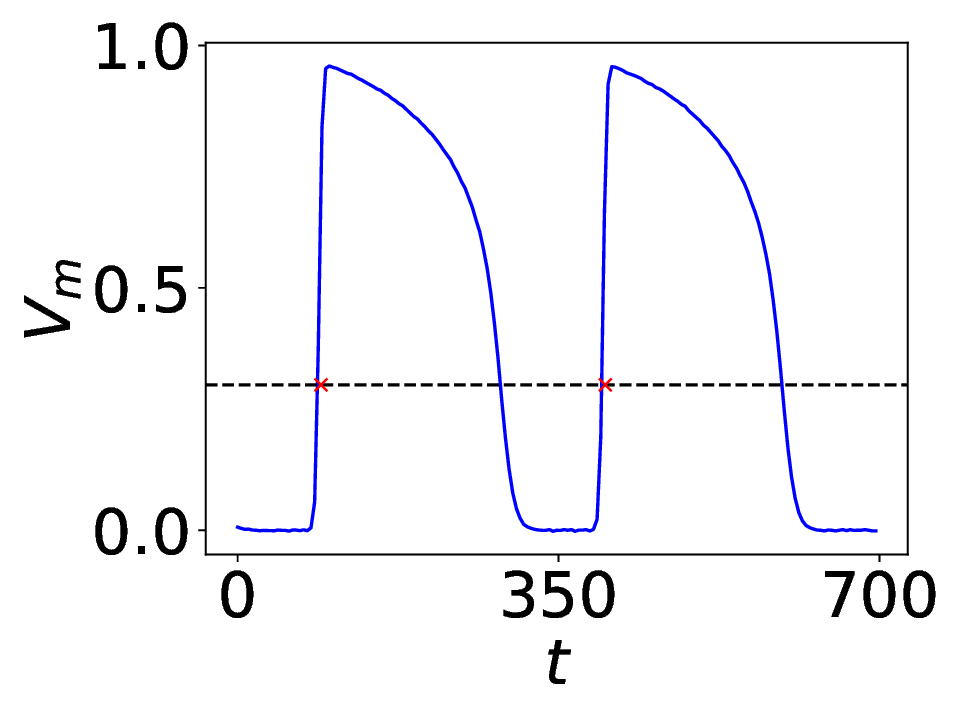}
    \includegraphics[width=0.49\linewidth]{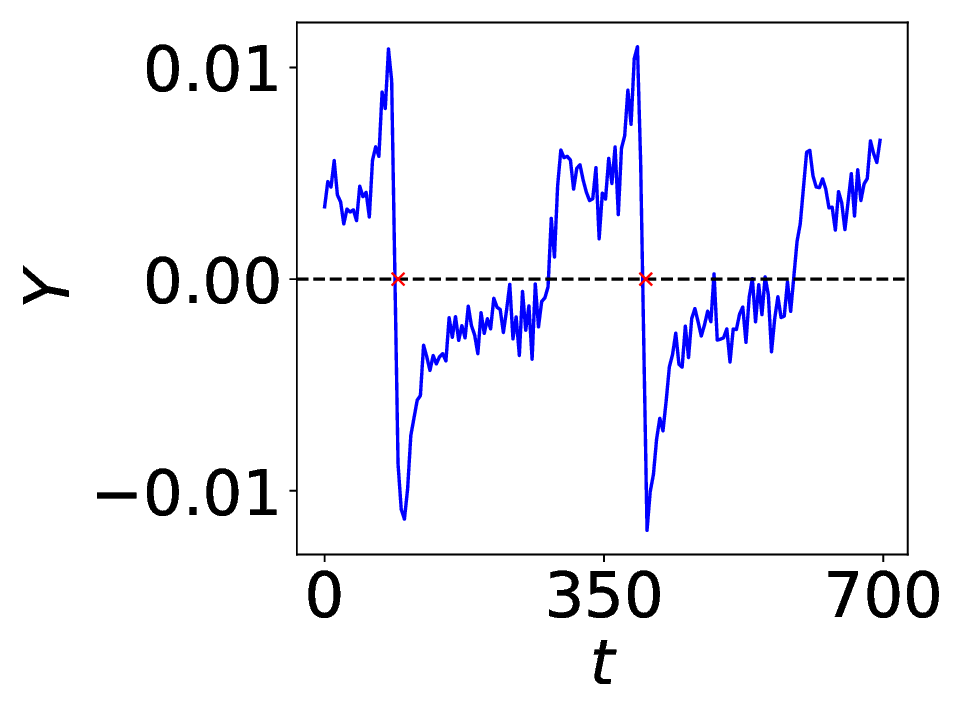}
\caption{Illustration of the calculation of the LAT. Left: the crossings (red crosses) with positive descent of the time trace $V_m(x_4,t)$ (in blue) with the horizontal line at $0.3$. Right: the crossings (red crosses) with negative descent of the time trace $y(r_4,t)$ with the horizontal line at 0.0. }
     \label{fig:LAT}
 \end{figure}
\end{itemize}

\printbibliography
\end{document}